\documentclass{elsarticle}
\usepackage{amsthm,enumerate, graphicx, amsmath,verbatim}
\usepackage[mathscr]{euscript}
\newtheorem{thm}{Theorem}[section]

\newtheorem{lem}[thm]{Lemma}

\newtheorem{cnj}[thm]{Conjecture}
\newtheorem{exm}[thm]{Example}
 \renewcommand\ell{l}

\theoremstyle{definition}

\newtheorem{defn}[thm]{Definition}
\DeclareMathOperator{\la}{la}

    \makeatletter
    \def\ps@pprintTitle{%
       \let\@oddhead\@empty
       \let\@evenhead\@empty
       \let\@oddfoot\@empty
       \let\@evenfoot\@oddfoot
    }
    \makeatother
    
\begin{document}
\title{Edge Coloring Signed Graphs}
\author{Richard Behr}
\address{Department of Mathematics, Binghamton University, Binghamton, NY 13902}
\ead{behr@math.binghamton.edu}

\begin{abstract} We define a method for edge coloring signed graphs and what it means for such a coloring to be proper. Our method has many desirable properties: it specializes to the usual notion of edge coloring when the signed graph is all-negative, it has a natural definition in terms of vertex coloring of a line graph, and the minimum number of colors required for a proper coloring of a signed simple graph is bounded above by $\Delta+1$ in parallel with Vizing's Theorem. In fact, Vizing's Theorem is a special case of the more difficult theorem concerning signed graphs. 
\end{abstract}
\maketitle

\section{Introduction}

A \textit{signed graph} is a graph in which each edge is assigned either a positive or negative sign. Signed graphs were invented by Harary in 1953 in order to help study a question in social psychology \cite{harary} and since then have proved to be a natural generalization of ordinary graphs in many ways. Typically the signed graph version of a graph theoretic structure is similar to the unsigned version, but more complicated due to the presence of the edge signs.  Various phenomena that are unseen in the world of ordinary graphs will manifest themselves in the signed world, leading to interesting insights into both ordinary graphs and signed graphs. For example, Zaslavsky \cite{zaslav2} discovered a theory of signed vertex coloring complete with deletion/contraction recurrence and chromatic polynomials that specializes to ordinary graphs when the signed graph is all positive. Zaslavsky's coloring construction makes use of ``signed colors", which are a new twist that must be introduced to define vertex coloring of a signed graph in an interesting fashion. In this paper we will use these same signed colors to construct a theory of signed edge coloring that has many nice properties, including compatibility with Zaslavsky's signed vertex coloring.

It is a well-known theorem of Vizing that the number of colors needed to properly edge color an ordinary simple graph is either $\Delta$ or $\Delta+1$ \cite{vz}. Thus there are two kinds of graphs---those that can be $\Delta$-colored (class $1$), and those that require one additional color (class $2$). In what follows, we will prove that every signed simple graph can be edge colored with $\Delta$ or $\Delta+1$ colors (Theorem \ref{goodvizing}), rendering ordinary Vizing's Theorem a special case of our new theorem concerning signed graphs. Interestingly, it is possible to change the class of a signed graph by modifying its signature, so that two signed graphs on the same underlying graph may fall into distinct sides of the dichotomy.

To prove Signed Vizing's Theorem we will employ the signed graph analogues of several techniques that are commonly used when studying Vizing's theorem, such as Kempe chains and the ``fan" used to recolor edges locally. Signed Kempe chains are particularly interesting, as they exhibit certain properties that are unseen in ordinary Kempe chains.

It is reasonable to expect nuances or exceptions to arise when generalizing an ordinary graph theory concept to signed graphs. For example, in ordinary graph theory the well-known Brooks' Theorem states that the number of colors needed to vertex color a connected graph is bounded above by $\Delta$ with two exceptions---complete graphs and odd cycles require one additional color. The signed version of Brooks' Theorem is nearly identical, but it turns out that negative cycles of even length, in addition to positive complete graphs and positive even cycles, are a third exception \cite{maca}. Remarkably, no additional exceptions arise in the signed version of Vizing's theorem, and the upper bound for the number of colors is $\Delta+1$ for both ordinary and signed graphs alike.

\section{Graphs and Signed Graphs}
\subsection{Graphs}

We write $\Gamma$ for a \emph{graph} and we write $V(\Gamma)$ and $E(\Gamma)$ for its vertex and edge sets respectively. Throughout, we will assume that every graph edge has two distinct endpoints, and that no two edges have the same pair of endpoints. In other words, we assume that all graphs are \emph{simple}. Often we write $e{:}vw$ for an edge $e$ with endpoints $v$ and $w$. If vertices $v$ and $w$ are connected by an edge we say that they are \emph{adjacent} or that they are \emph{neighbors}.

An \emph{incidence} of $\Gamma$ is a pair $(v,e)$ such that vertex $v$ is an endpoint of edge $e$. The set of all incidences of $\Gamma$ is written $I(\Gamma)$. If we write $(v,vw)$ it is understood that we are referring to the incidence between $v$ and edge $e{:}vw$.

A \emph{circle} is a connected $2$-regular subgraph. A \emph{path} is a sequence of adjacent vertices and connecting edges that never repeats an edge or a vertex. A \emph{trail} has the same definition as a path, except that a trail may repeat vertices (but not edges). Thus, every path is a trail, but not every trail is a path. For trails and paths, we call $v_0$ and $v_n$ the \emph{endpoints}, while the other vertices are \emph{interior vertices}. In a path the endpoints and interior vertices are distinct, but in a trail there may be interior vertices that are also endpoints. Often, we specify a trail by listing its vertices in order inside of parenthesis. For a trail $T$ with endpoints $t_0$ and $t_n$ we write $T=(t_0,...,t_n)$.

A \emph{matching} $M$ in $\Gamma$ is a collection of edges of $\Gamma$ such that no two edges of $M$ share an endpoint. An \emph{independent set} $J$ in $\Gamma$ is a collection of vertices such that no edge has both endpoints in $J$.

\subsection{Signed Graphs}

A \emph{signed graph} is a pair $\Sigma = (\Gamma,\sigma)$, where $\Gamma$ is a graph and $\sigma : E(\Gamma) \rightarrow \{+,-\}$ is the \emph{signature}. We write $|\Sigma| = \Gamma$ for the \emph{underlying graph} of $\Sigma$, the unsigned graph obtained by forgetting all of the signs. 

A circle in $\Sigma$ is \emph{positive} if the product of its edge signs is positive, and \emph{negative} otherwise. A subgraph of $\Sigma$ is \emph{balanced} if each of its circles is positive, and \emph{unbalanced} otherwise. A subgraph of $\Sigma$ is \emph{antibalanced} if each of its circles is either positive and even in length or negative and odd in length.

We write $-\Sigma$ for the \emph{negation} of $\Sigma$, defined by $-\Sigma = (|\Sigma|, -\sigma)$. It is easy (even easier, once we define switching) to check that $\Sigma$ is balanced if and only if $-\Sigma$ is antibalanced. 

\emph{Switching} $\Sigma$ by $v \in V(\Sigma)$ means negating the sign of every edge that has $v$ as an endpoint. Switching $\Sigma$ by $X \subseteq V(\Sigma)$ means switching each $v \in X$ in turn. If $\Sigma'$ is obtained from $\Sigma$ by switching, we say they are \emph{switching equivalent}, written $\Sigma' \sim \Sigma$. It is straightforward to check that switching equivalence is indeed an equivalence relation. The \emph{switching class} of $\Sigma$ is the equivalence class of $\Sigma$ under this equivalence relation and is denoted by $[\Sigma]$. 

A fundamental theorem concerning switching (as found in \cite{zaslav3}) is that $\Sigma' \sim \Sigma$ if and only if $|\Sigma'|=|\Sigma|$ and $\Sigma'$ and $\Sigma$ have the same balanced circles. Thus, any property of signed graphs that depends only on the signs of the circles is invariant for all graphs contained in $[\Sigma]$. As we will see, signed edge coloring using a certain number of colors is an example of this phenomenon. 

\section{Edge Colorings}

In this section we will give a natural definition for edge coloring a signed graph. Recall that an \emph{edge coloring} of an ordinary graph $\Gamma$ is an assignment of colors (typically elements of $\{1,\ldots,n\}$) to its edges. Such a coloring is \emph{proper} if no two adjacent edges receive the same color. Our definition is similar, but we define edge coloring in terms of incidences (rather than edges themselves) in order to incorporate edge signs.

To edge color a signed graph, we need a more sophisticated set of colors than $\{1,\ldots,n\}$. Let $M_n=\{0,\pm 1,\ldots,\pm k\}$ if $n=2k+1$, and $M_n=\{\pm 1,\ldots,\pm k\}$ if $n=2k$. The $M_n$ are called \emph{signed color sets} and they contain \emph{signed colors}. The colors $+a$ and $-a$ have the same \emph{magnitude}, but are \emph{opposite}. These are the same signed color sets used in both \cite{zaslav2} and \cite{maca} to study signed vertex coloring.

\begin{figure}[h!]
\centering
\includegraphics[scale=.07]{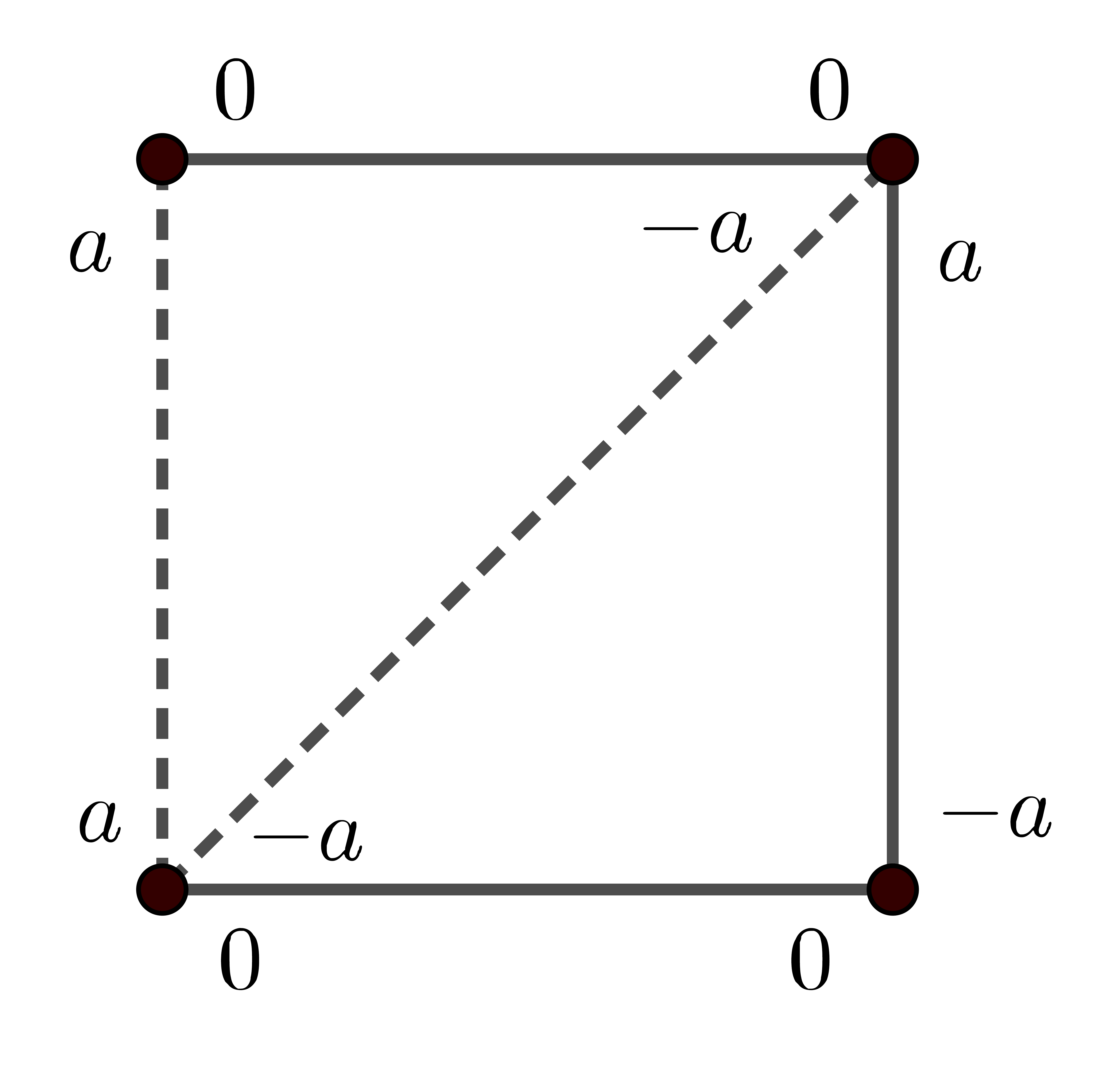}
\caption[A signed graph and a proper $3$-coloring.]{Here we see a signed graph and a proper $3$-coloring. Throughout, we indicate negative edges with dashed lines and positive edges with solid lines.}
\label{fig1}
\end{figure}

\begin{defn}An \emph{$n$-edge coloring} (or more briefly, an \emph{$n$-coloring}) $\gamma$ of $\Sigma$ is an assignment of colors from $M_n$ to each vertex-edge incidence of $\Sigma$ subject to the condition that $\gamma(v,e) = -\sigma(e) \gamma(w,e)$ for each edge $e{:}vw$.

An $n$-coloring is \emph{proper} if for any two incidences $(v,e)$ and $(v,f)$ involving the same vertex, $\gamma(v,e) \neq \gamma(v,f)$.
\end{defn}
Intuitively, negative edges act like unsigned edges since they have the same color at both incidences. However, positive edges act differently and instead have opposite colors at their incidences. See Figure \ref{fig1} for an example of a proper coloring. 

Since both incidences of a negative edge $e{:}vw$ receive the same color, we will sometimes write $\gamma(e) = \gamma(v,e) = \gamma(w,e)$. When such notation is used it is understood that the edge in question is negatively signed. We cannot afford this luxury when it comes to positive edges.

The color $0$ is its own opposite, so when the color $0$ is available we are allowed to color both positive and negative edges with $0$ at both of their incidences. The subsequent theory often becomes simpler when the color $0$ is not available, so it is sometimes convenient to discuss \emph{zero-free} colorings---those which omit $0$ and hence have an even number of colors available.

\begin{defn}Let $\gamma$ be an $n$-coloring of $\Sigma$. If there exists an edge $e$ such that $\gamma(v,e) = a$, then we say the color $a$ is \emph{present} at $v$. Otherwise, $a$ is \emph{absent} at $v$. 
\end{defn}

To rephrase the definition of a proper $n$-coloring, $\gamma$ is proper if and only if each color from $M_n$ is present at each vertex at most once. 

\subsection{Basic Coloring Properties}

We begin by reiterating that negative edges receive the same color at both of their incidences. Hence, colored negative edges are essentially the same as colored unsigned edges. We use this to record the following lemma.

\begin{lem}\label{specialize} If every edge of $\Sigma$ is negative, then there is a one-to-one correspondence between (proper) $n$-colorings of $\Sigma$ and (proper) $n$-colorings (in the usual unsigned sense) of $|\Sigma|$.
\end{lem}

The fact that signed edge coloring specializes to ordinary edge coloring when $\Sigma$ is all negative (as opposed to all positive, as one might expect) is a consequence of the definition of the signed line graph, as we will see in later sections. 

Another basic feature of signed edge coloring is its compatibility with switching. 

\begin{lem}\label{switchin} Suppose $\gamma$ is a proper $n$-coloring of $\Sigma=(\Gamma,\sigma)$ and suppose $\Sigma'=(\Gamma,\sigma')$ is obtained from $\Sigma$ by switching a vertex set $X$. Define a new coloring $\gamma \,'$ which is obtained from $\gamma$ by negating all colors on all incidences involving vertices from $X$. Then, $\gamma \,'$ is a proper $n$-coloring of $\Sigma'$.

\begin{proof} We will describe the case where we switch a single vertex, since all larger cases can be considered one vertex at a time. Consider a vertex $v \in X$. Since $\gamma$ is an edge coloring, $\gamma(v,e)=-\sigma(e)\gamma(w,e)$ for any edge $e{:}vw$ incident with $v$. Switching $v$ changes the signs of all edges adjacent to $v$. Hence, we have $\gamma \,'(v,e)=-\sigma'(e)\gamma \,'(w,e)$ and thus $\gamma \,'$ is an edge coloring. We further see that $\gamma\, '$ is proper, since $\gamma(v,e)\neq \gamma(v,f)$ implies that $\gamma \,'(v,e) \neq \gamma \,'(v,f)$. 
\end{proof}
\end{lem}

\begin{figure}[h!]
\centering
\includegraphics[scale=.05]{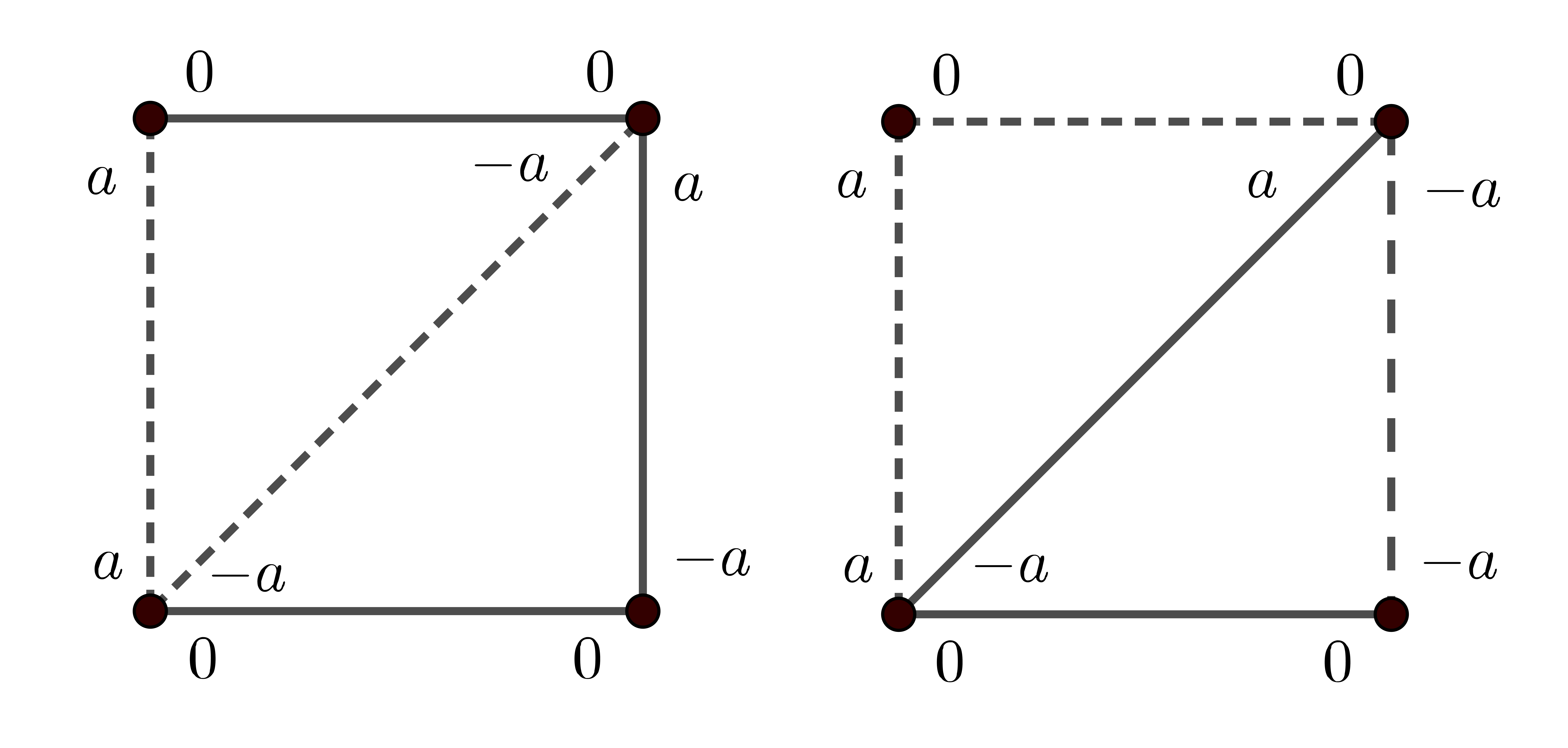}
\caption[Applying a switching function to an edge colored signed graph.]{On the left-hand graph, we have the same proper coloring as seen in Figure 1. On the right-hand graph, we have switched the top-right vertex and negated all colors incident with this vertex to obtain a proper coloring of the switched graph.}
\label{fig2}
\end{figure}

Thus, if we obtain a proper $n$-coloring of $\Sigma$, we automatically obtain a proper $n$-coloring for every member of $[\Sigma]$. Lemma \ref{switchin} is illustrated in Figure \ref{fig2}.  

\subsection{Magnitude Subgraphs}

In this section we study the properties of subgraphs whose edges are colored with a single magnitude. Let $\Sigma$ be a signed graph and let $\gamma$ be a proper $n$-coloring of $\Sigma$. We write $\Sigma_{a}[\gamma]$ for the set of edges of $\Sigma$ that are colored using $\pm a$ with respect to $\gamma$. If there is only one coloring that we have in mind, we write $\Sigma_a$. We call $\Sigma_a$ the \emph{$a$-graph} of $\Sigma$ with respect to $\gamma$.

We first observe that $\Sigma_a$ has maximum degree $2$, since at most $a$ and $-a$ are present at each vertex of $\Sigma_a$. Thus, each component of $\Sigma_a$ is either a path or a circle. When $a=0$ the maximum degree is $1$, and hence $\Sigma_0$ is a matching. We now describe which kinds of paths and cycles can possibly appear in $\Sigma_a$ when $a\neq 0$.

\begin{lem} Every signed path can be properly edge colored with $\pm a$ (where $a \neq 0$). Furthermore, every signed path has exactly two different $\pm a$ colorings.
\begin{proof} Switch so that the path is all negative. Then, color the edges of the path so that they alternate between $-a$ and $a$. Finally, switch back to the original signature of the path, negating colors at the switched vertices as in Lemma \ref{switchin}. Clearly there are exactly two possibilities for each path---simply negate all colors to change between the two possible $\pm a$ colorings.
\end{proof}
\end{lem}

\begin{figure}[h!]
\centering
\includegraphics[scale=.069]{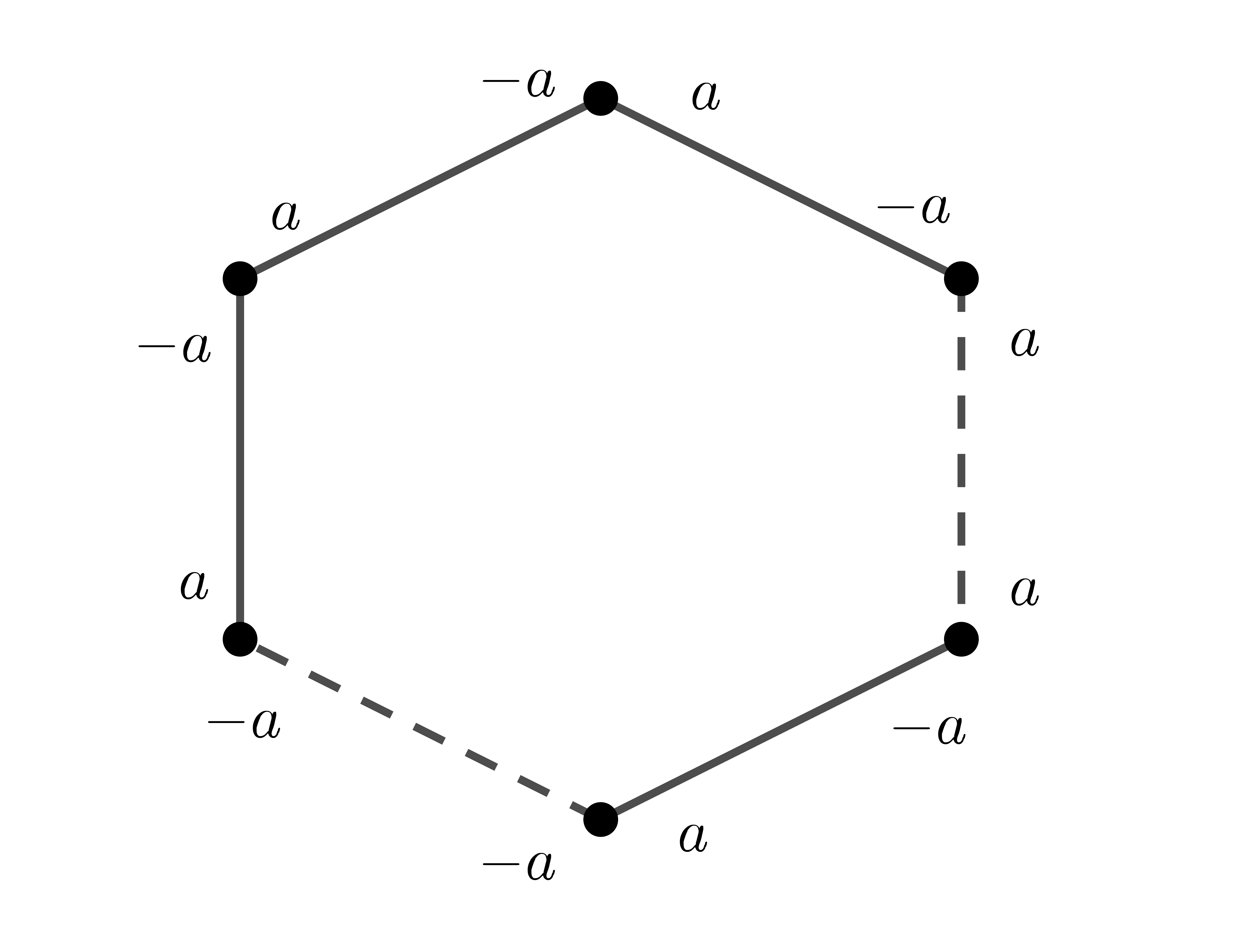}
\caption{A positive circle and one of its two possible $2$-colorings. \label{fig3}}
\end{figure}

\begin{lem}A signed circle $C$ can be properly colored with $\pm a$ ($a \neq 0$) if and only if $C$ is positive. Furthermore, every positive circle has exactly $two$ $\pm a$ colorings.
\begin{proof} First, suppose $C$ can be properly colored with $\pm a$. Choose an edge $e{:}vw$ of $C$ and switch so that $C{\setminus}e$ is an all-positive path from $v$ to $w$. Since $C{\setminus}e$ is all-positive, up to choice of names $a$ is present at $v$ in $C{\setminus}e$ and $-a$ is present at $w$ in $C{\setminus}e$. Hence, $e$ must be positive, or else a coloring using $\pm a$ would be impossible.

Conversely, any positive circle can be colored with $\pm a$. Simply switch so that the circle is all positive and color incidences alternating between $a$ and $-a$ around the circle. 

Finally, we note that given any $2$-coloring of a positive circle, we can obtain another $2$-coloring by negating all of the colors. Clearly, every positive circle has just two possible $2$-colorings. A $2$-colored positive circle is shown in Figure \ref{fig3}.
\end{proof}
\end{lem}

Therefore in a proper edge coloring, $\Sigma_a$ consists of paths and positive circles so that $\Sigma_a$ is balanced. Thus, a proper coloring corresponds to a partition of the edges of $\Sigma$ into balanced subgraphs of maximum degree $2$ and a single matching (if the color $0$ is available). This is not a one-to-one correspondence, as each connected component in such a partition may be colored in one of two ways (except for edges that are colored $0$). 

We see the phenomenon mentioned above in ordinary graphs as well. Suppose $c$ is an ordinary edge coloring of $\Gamma$ (using $M_n$ for the color set). In this case, $\Sigma_a$ ($a \neq 0$) is a bipartite subgraph of maximum degree $2$ (or a matching if $a=0$). 

When $\Sigma$ is all negative, the balanced subgraphs are precisely the bipartite subgraphs. Hence, the partition into bipartite subgraphs induced by an ordinary edge coloring is a special case of the balanced partition induced by a signed edge coloring.

\section{Vizing's Theorem for Signed Graphs}
For an ordinary graph $\Gamma$, we write $\chi'(\Gamma)$ for the \emph{chromatic index} of $\Gamma$---the minimum number of colors used in any proper edge coloring of $\Gamma$. It is a classic result of Vizing that every simple ordinary graph $\Gamma$ satisfies $\Delta(\Gamma) \leq \chi'(\Gamma) \leq \Delta(\Gamma)+1$, where $\Delta(\Gamma)$ is the maximum degree of a vertex in $\Gamma$ \cite{vz}.

For a signed graph $\Sigma$ we borrow the existing notation and write $\chi'(\Sigma)$ for its \emph{chromatic index}, which we define to be the smallest $n$ such that $\Sigma$ has a proper edge coloring using colors from $M_n$. As an easy application of Lemma \ref{specialize}, we see that $\chi'(|\Sigma|)=\chi'(\Sigma)$ if $\Sigma$ is all negative. Because of this, our use of the usual notation is justified.

The lower bound for $\chi'(\Sigma)$ is $\Delta(\Sigma)$ since there are $\Delta(\Sigma)$ different incidences to be colored at a vertex of maximum degree. Our goal is now to prove an upper bound of $\Delta(\Sigma)+1$, but we first require some new machinery. We will define a signed version of a Kempe chain, which is a tool that is frequently used in the study of ordinary edge coloring.

\subsection{Signed Kempe Chains}

The main difference between ordinary and signed Kempe chains is that the signed version is a trail (which may intersect itself at a vertex), and not a path.

We begin with a convenient notational definition. If $T = (v_0,\ldots,v_m)$ is a signed trail we write $t_k$ for the number of positive edges that appear on $T$ between $v_0$ and $v_k$. 

\begin{defn}\label{kempe}Suppose $\gamma$ is a proper $n$-coloring of $\Sigma$. If $a$ is absent at a vertex $v_0$ and $b$ is present at $v_0$, we define the \emph{$a/b$-chain at $v_0$} to be the maximal trail $T=(v_0,\ldots,v_m)$ starting at $v_0$ with the properties:
\begin{enumerate}
\item The edge magnitudes alternate between $|a|$ and $|b|$ along $T$ (starting with $|b|$).
\item $\{\gamma(v_i, v_{i-1}v_i), \gamma(v_i, v_iv_{i+1})\} = \{ (-1)^{t_i} a, (-1)^{t_i} b \}$ for all $i \neq 0,m$. 
\end{enumerate}
\end{defn}

We write $K_{a,b}(v_0, v_m)$ for the $a/b$-chain at $v_0$. In the following lemma we note that signed Kempe chains specialize to ordinary Kempe chains when $\Sigma$ is all negative.

\begin{lem}If $\Sigma$ is all negative, $K_{a,b}(v_0, v_m)$ is a path.

\begin{proof} Since all edges are negative, the second property in the definition of the signed Kempe chain becomes $\{\gamma(v_i, v_{i-1}v_i), \gamma(v_i, v_iv_{i+1})\} = \{a , b \}$ for all $i \neq 0,m$. Thus, the signed Kempe chain has maximum degree $2$ and hence does not intersect itself. It must be a path.
\end{proof}
\end{lem}

\begin{figure}[h!]
\centering
\includegraphics[scale=.05]{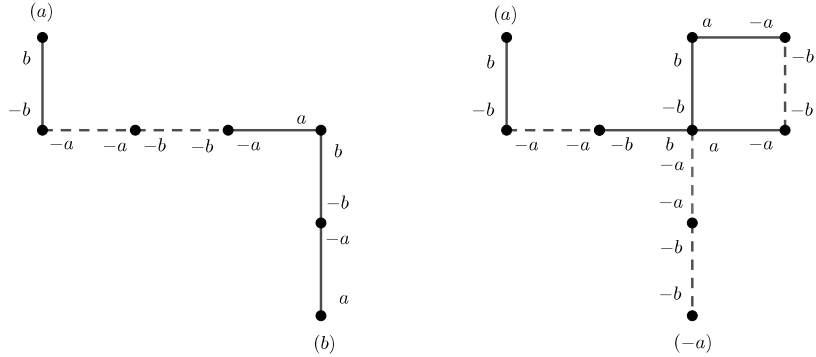}
\caption[Two examples of $a/b$-chains.]{On the left, an $a/b$-chain that does not intersect itself. On the right, an $a/b$-chain that does intersect itself. Parenthetical colors indicate their absence at the specified vertex.}
\label{fig4}
\end{figure}

While an ordinary Kempe chain does not intersect itself, a signed Kempe chain may intersect itself. For example, let $(v_0,\ldots,v_m)$ be the signed $a/b$-chain at $v_0$. Suppose at vertex $v_j$ we have $\gamma(v_j, v_{j-1}v_j) = a$ and $\gamma(v_j, v_jv_{j+1}) = b$. Then, vertex $v_j$ may appear once again in the Kempe chain (say, $v_j = v_s$, $s > j$), as long as $\gamma(v_s, v_{s-1}v_s) = -a$ and $\gamma(v_s, v_s v_{s+1}) = -b$ (or vice versa). There is at most one self intersection at each vertex, since the only available colors are $\pm a$ and $\pm b$. A self-intersecting Kempe chain is shown in Figure \ref{fig4}. The following lemma shows that in the event of a self-intersection, the path between $v_j$ and $v_s$ must have a certain sign.

\begin{lem} Suppose $K_{a,b}(v_0, v_m)$ has a vertex $v_j = v_s$ that appears twice $(s > j)$. Let $Q$ be the subtrail $(v_j,\ldots,v_s)$. Then, $Q$ has an odd number of positive edges.
\begin{proof} If $Q$ has an even number of positive edges, $\{\gamma(v_j, v_{j-1}v_j), \gamma(v_j, v_jv_{j+1})\} = \{\gamma(v_s, v_{s-1}v_s), \gamma(v_s, v_sv_{s+1})\}$, which contradicts the definition of a proper coloring.
\end{proof}
\end{lem}

The main purpose of signed Kempe chains is that we can use them to change the colors present at a given vertex.

When $a,b\neq0$, we modify $K_{a,b}(v_0, v_m)$ by performing the $a/b$-\emph{swap} at $v_0$, which is the act of interchanging $a$ with $b$ and $-a$ with $-b$ for all incidences in $K_{a,b}(v_0, v_m)$. Thus, the $a/b$ swap changes the color present at $v_0$ from $b$ to $a$, and also changes the color present at $v_m$.

Intuitively, we imagine that changing $\gamma(v_0, v_0v_1)$ from $b$ to $a$ creates a sequence of subsequent changes that must be made in order to preserve the propriety of the coloring. This sequence of changes propagates along $K_{a,b}(v_0,v_m)$ and ends when $K_{a,b}(v_0,v_m)$ ends, allowing us to maintain a proper coloring. We now describe in detail what happens when performing a swap.

\begin{lem}\label{swap1}Consider $K_{a,b}(v_0, v_m)$ where $a \neq 0$ and $b \neq 0$. Performing the $a/b$-swap at $v_0$ does not change the present or absent colors at any vertex except for $v_0$ and $v_m$. It changes the present colors at $v_0$ and $v_m$, interchanging $a$ with $b$ at $v_0$ and $(-1)^{t_m}a$ with $(-1)^{t_m}b$ at $v_m$.
\end{lem}

We now discuss what happens when a signed Kempe chain involving the color $0$ intersects itself. Since $+0 = -0$, such a chain has maximum degree $3$ and hence does not behave the same as a zero-free chain.

\begin{lem} Suppose $K_{a,0}(v_0,v_m)$ has its first self-intersection at $v_j=v_s$. Then, $v_s = v_m$. In other words, the Kempe chain must terminate at its first self-intersection.
\begin{proof} Suppose without loss of generality that $\gamma(v_j, v_{j-1}v_j) = a$ and $\gamma(v_j,v_jv_{j+1}) = 0$. Then the only possibility for $\gamma(v_s, v_{s-1}v_s)$ is $-a$. However, the chain cannot continue past $v_s$ since all three available colors are already present at $v_s$. 
\end{proof}
\end{lem}

Finally, we make note of what happens if we swap a self-intersecting Kempe chain that involves the color $0$. 

\begin{figure}[h!]
\centering
\includegraphics[scale=.081]{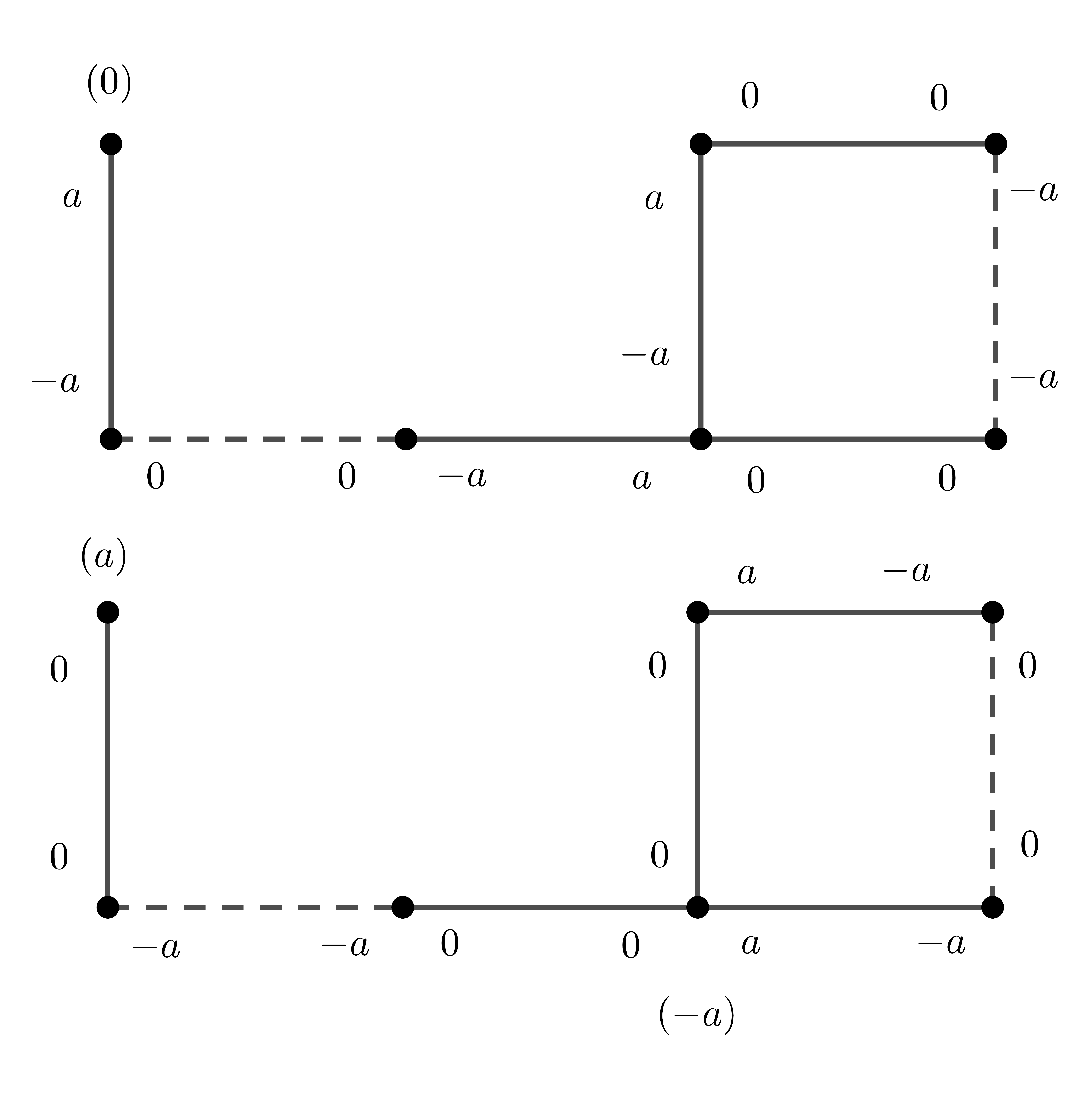}
\caption[Swapping an $a/0$-chain.]{On the top, an $a/0$-chain that terminates when it intersects itself. On the bottom, the result of swapping the chain. Swapping the chain ruins the propriety of the coloring.}
\label{fig5}
\end{figure}

\begin{lem}Consider $K_{a,0}(v_0,v_m)$ where $a \neq 0$. If $K_{a,0}(v_0,v_m)$ intersects itself at $v_j = v_m$ where $j < m$, then performing the $a/0$-swap at $v_0$ creates a single impropriety at $v_m$. Namely, either $\gamma(v_m, v_{m-1}v_m)=\gamma(v_j, v_{j-1}v_j) = 0$ or $\gamma(v_m, v_{m-1}v_m)=\gamma(v_j, v_{j}v_{j+1}) = 0$. 
\begin{proof} Since $a, -a$, and $0$ are present at $v_m$ before the swap, then $0$, $0$, and either $a$ or $-a$ are present at $v_m$ after the swap. No two incidences involving the same vertex may be colored $0$ in a proper coloring.
\end{proof}
\end{lem}
Because swapping a self-intersecting chain involving $0$ creates an improper coloring, we will consider only zero-free Kempe chains in subsequent arguments and deal with $0$ a different way. The above lemma is illustrated in Figure \ref{fig5}.
\subsection{Signed Vizing's Theorem}
We are now ready to prove the signed generalization of Vizing's theorem. We will first prove the zero-free version of the theorem (Theorem \ref{vizing}) using a method that involves Kempe chains. The zero-free version gives the desired upper bound of $\Delta(\Sigma)+1$ when $\Delta(\Sigma)$ is odd, and a weaker upper bound of $\Delta(\Sigma)+2$ when $\Delta(\Sigma)$ is even. This occurs because it is impossible to use $\Delta(\Sigma)+1$ colors when $\Delta(\Sigma)$ is even and the color $0$ is unavailable. After proving Theorem \ref{vizing}, we incorporate the color $0$ using a method that does not involve Kempe chains, bringing the upper bound down to $\Delta(\Sigma)+1$ when $\Delta(\Sigma)$ is even.

The statement of our main theorem is as follows.

\begin{thm}[Signed Vizing's Theorem]\label{goodvizing} For a signed simple graph $\Sigma$, $\Delta(\Sigma) \leq \chi'(\Sigma) \leq \Delta(\Sigma)+1$.
\end{thm}

We first define a \emph{fan}, which is a device that allows us to manipulate the colors locally at a given vertex without affecting the propriety of a coloring. 

\begin{defn} Let $\Sigma$ be a signed graph and let $e{:}uv_0$ (called the \emph{initial edge}) be an edge of $\Sigma$. Let $\gamma_0$ be a proper edge coloring (called the \emph{initial coloring}) of $\Sigma {\setminus} e$. Assume by switching that all edges incident with $u$ are negative.  Let $v=v_0,\ldots,v_s$ be a maximal sequence of neighbors of $u$ such that $\gamma_0(v_i, uv_i)$ is absent at $v_{i-1}$ for all $1 \leq i \leq s$. The \emph{fan} at $u$ is the subgraph induced by all of the $uv_i$ edges. The edges of the fan are written $e_i:=uv_i$. We say that $u$ is the \emph{hinge} of the fan. 
\end{defn}
The purpose of a fan is to allow us to interchange the colors on the edges adjacent to $u$. The colorings we obtain by interchanging colors are called \emph{shifted colorings} and are defined as follows.
\begin{defn} Let $F$ be a fan with initial coloring $\gamma_0$ and edges $e_0,\ldots,e_s$ (where $e_0{:}uv_0$ is the uncolored initial edge). We define a sequence of \emph{shifted colorings}, $\gamma=\gamma_1,\ldots,\gamma_s$, such that:

\begin{enumerate}
\item The edge $e_i$ is not colored in $\gamma_i$.
\item $\gamma_i(u,e_j)=\gamma_i(v_j,e_j) = \gamma_0(v_{j+1},e_{j+1})$ for $j \in \{0,\ldots,i-1\}$.
\item $\gamma_i=\gamma_0$ otherwise.
\end{enumerate}  
\end{defn}

The reader is encouraged to envision $v_0$ as being at the bottom, $v_s$ as being at the top, and the $\gamma_i$ as being obtained by shifting the colors of the edges $e_1,\ldots,e_i$ downwards, leaving $e_i$ uncolored. By design, each of the $\gamma_i$ is proper and uses the same set of colors as $\gamma_0$. Also notice that the colors present or absent at $u$ are exactly the same in all of the $\gamma_i$. 

Under certain conditions, we will be able to use a fan in conjunction with Kempe chains to extend a proper coloring of $\Sigma{\setminus}e_0$ to a proper coloring of $\Sigma$ using the same set of colors. The following lemma forms the bulk of the proof of Theorem \ref{vizing}, and hence is a substantial portion of the proof of Theorem \ref{goodvizing}.

\begin{lem}\label{fan} Let $\Sigma$ be a signed graph and let $e_0{:}uv_0$ be one of its edges. Let $\gamma_0$ be a proper zero-free edge coloring of $\Sigma{\setminus}e_0$ using $n$ colors, and suppose that there is at least one color absent at $u$ and at each neighbor of $u$. Furthermore, suppose colors of the same magnitude are absent at $u$ and $v_0$. Then, there exists an $n$-coloring of $\Sigma$.
\begin{proof}
First, assume by switching that all edges incident with $u$ are negative. Let $a$ be a color absent at $u$, and suppose a color of the same magnitude is absent at $v_0$. If $a$ is also absent at $v_0$, we extend the coloring by setting $\gamma_0(e_0)=a$ and we are done. Otherwise, assume $-a$ is absent at $v_0$. If $-a$ is absent at $u$, then we can also extend the coloring, and so we assume $-a$ is present at $u$. So, both $u$ and $v_0$ have degree $1$ in $(\Sigma{\setminus}e)_a$.

Now we build a fan $F$ with hinge $u$ and initial edge $e_0$, and with shifted colorings $\gamma_0,\ldots,\gamma_s$. Let $e_0,\ldots,e_s$ be the edges of $F$. Since $-a$ is present at $u$ and absent at $v_0$, we choose $e_1$ to be the edge adjacent to $u$ with color $-a$. 

Next, let $b$ be a signed color that is absent at $v_s$ with respect to $\gamma_0$. Then, $b$ is absent at $v_s$ with respect to all of the $\gamma_i$.  If $b$ is absent at $u$, then we can simply extend $\gamma_s$ by setting $\gamma_s(e_s) = b$. So, we assume that $b$ is present at $u$. Since $F$ is maximal (by definition), there must be some $1 \leq j \leq s-1$ such that $\gamma_0(e_j) = b$. See Figure \ref{fig6} for an illustration of $F$.

\begin{figure}[h!]
\centering
\includegraphics[scale=.050]{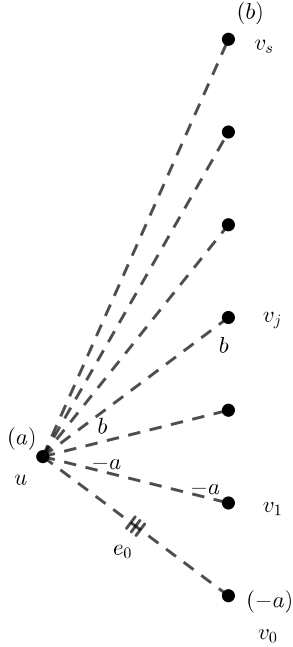}
\caption[A depiction of $\gamma_0$ and $F$.]{A depiction of $\gamma_0$ and $F$. The triple marking on the initial edge $e_0$ indicates that it is currently uncolored.\label{fig6}}
\end{figure} 

Now we consider the final shifted coloring $\gamma_s$. If $a$ is absent at $v_s$ with respect to $\gamma_s$, we simply extend $\gamma_s$ by coloring $\gamma_s(e_s)=a$. Thus, we assume that $a$ is present at $v_s$. Since $b$ is absent at $v_s$, there is an $a/b$-chain starting at $v_s$, which we denote by $T$. We write the vertices of $T$, in order, as $v_s=t_0,t_1,\ldots,t_r$, so that $T=(v_s,t_1,...,t_r)$. If we swap the $a/b$-chain at $v_s$ we can color both incidences of $e_s$ with the color $a$, unless performing the $a/b$-swap at $v_s$ makes it so that $a$ is present at $u$. By Lemma \ref{swap1}, this can only happen if $u=t_r$ (i.e., $T$ ends at $u$), which we now assume. So, to use our Kempe chain notation, $T = K_{a,b}(v_s,u)$ with respect to $\gamma_s$. It is worth noting that while $T$ ends at $u$, it is also possible for $T$ to pass through $u$ once before ending at $u$. If $T$ does pass through $u$ before ending at $u$, it must pass through $u$ at consecutive edges whose colors are $-a$ and $-b$ (or vice versa). Thus, there may be $1$, $2$, or $3$ edges of $T$ that are contained in $F$.

Once again consider $\gamma_s$. Since performing the $a/b$-swap at $v_s$ makes the color $a$ present at $u$, the last incidence of $T$ must be colored $b$ . Thus, the last edge of $T$ is $uv_{j-1}$, which is the edge of $F$ that is colored $b$ with respect to $\gamma_s$. We now break the proof into several cases, depending on the nature of the edges in the intersection of $T$ and $F$. Let $X$ be the set of edges in the intersection of $T$ and $F$.

Case 1: The edge $e_{j-1}$ (colored $b$) is above all other edges of $X$ with respect to $\gamma_s$. By above, we mean that the index $j-1$ is greater than the index of all other edges in $X$.  Shift from $\gamma_s$ to $\gamma_{j-1}$, so that $e_{j-1}$ is uncolored and $b$ is absent at $v_{j-1}$. The key is that performing this shift only disturbs $T$ at its last edge $e_{j-1}$, since all other edges in $X$ are below $e_{j-1}$ and hence do not have their colors changed when shifting from $\gamma_s$ to $\gamma_{j-1}$. Now, we perform the $a/b$-swap at $v_{j-1}$, which propagates backwards along the trail $(v_{j-1}=t_{r-1},t_{r-2},\ldots,t_1,t_0=v_s)$ and terminates at $v_s$. This changes the color present at $v_s$ from $a$ to $b$ and also changes the color absent at $v_{j-1}$ from $b$ to $a$ without changing any other present or absent colors (by Lemma \ref{swap1}). Thus, we are now free to color both incidences of $e_{j-1}$ with the color $a$, completing this case. 

We pause to note two things. First, if $b=-a$ then the only edge in $X$ is $e_{j-1}$. Thus we have disposed of the case where $b=-a$  in the previous paragraph. In what follows, we will assume $b \neq -a$. Second, the cases where $|X|=1$ and $|X|=2$ are proved in the previous paragraph. Indeed, if $|X|=1$ then $T$ must end at $e_{j-1}$, and if $|X|=2$ then $X$ must contain only $e_0$ (which is colored $-a$ with respect to $v_s$) and $e_{j-1}$. To see this, suppose that $|X|=2$ and the edges of $X$ are $e_k$ (colored $-b$) and $e_{j-1}$ (colored $b$). Then either $T=(t_0,\ldots, v_k, u,\ldots u)$ or $T=(t_0,\ldots, u,v_k,\ldots u)$. In the first case, the incidence after $(u, uv_k)$ along $T$ must be colored $-a$, since $(u,uv_k)$ is colored $-b$. In the second case, the incidence before $(u, uv_k)$ must be colored $-a$. Thus, in both cases the edge $e_0$ (which is colored $-a$ and is contained in $F$) must be contained in $T$. This contradicts the fact that $|X|=2$. Thus, if $|X|=2$ it must contain $e_{j-1}$ and $e_0$ and hence Case 1 applies.

We point out that if $\Sigma$ is all negative then our work so far essentially implies ordinary Vizing's Theorem. One can prove ordinary Vizing's Theorem using the techniques described above, but the only possible case is where $|X|=1$. Case 1 is illustrated in Figure \ref{fig7}.
\begin{figure}[h!]
\centering
\includegraphics[scale=.056]{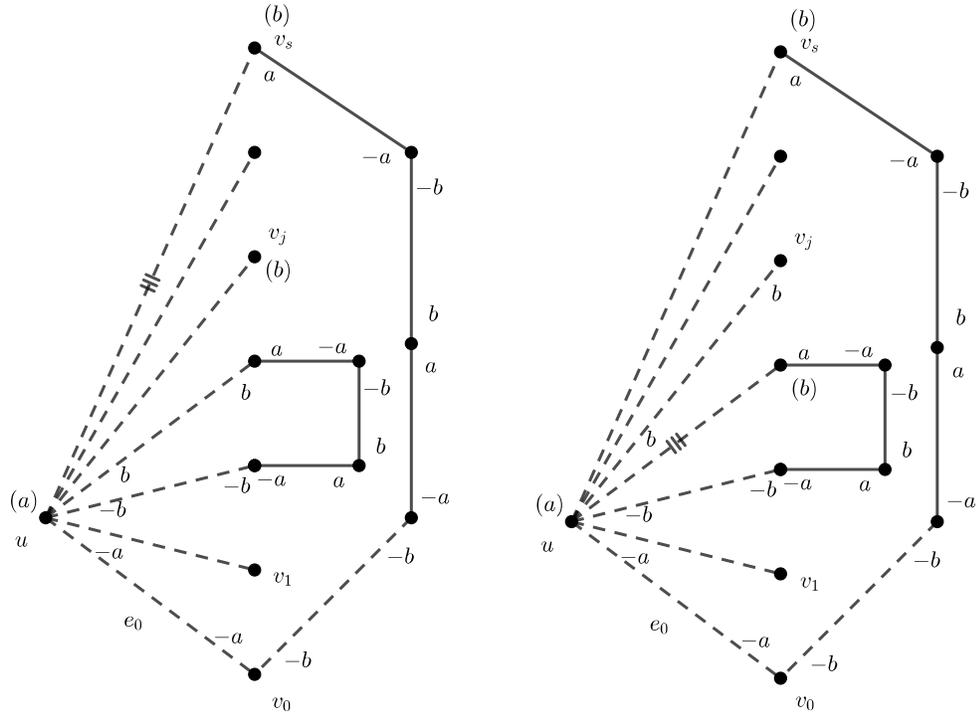}
\caption[An illustration of Case 1 in Theorem \ref{goodvizing}.]{An illustration of Case 1. On the left, we see $\gamma_s$ along with $T$, which ends at $u$ and intersects the fan at three edges. The $b$-colored edge is above the other two. On the right, we have shifted up to $\gamma_{j-1}$. Performing this shift does not break the Kempe chain exept at its final edge (the $b$-colored edge). On the right-hand graph, performing the $a/b$-swap at $v_{j-1}$ allows us to color $e_{j-1}$ with color $a$. Once again, the triple marking indicates that the given edge is uncolored.}
\label{fig7}
\end{figure} 

Case 2: The edge $e_{j-1}$ is not above all other edges of $X$ with respect to $\gamma_s$. Due to the previous discussion, this implies that $|X|=3$. Let the edges in $X$ be $e_k$ (colored $-b$), $e_0$ (colored $-a$), and $e_{j-1}$ (colored $b$). Since $\gamma_s(e_0)=-a$, we see that $e_k$ is above all the other edges. In other words, $k > j-1 >0$.

 To prove Case 1 above, we relied on the fact that shifting to $\gamma_{j-1}$ disturbed $T$ only at its final edge. In Case 2 we can no longer shift to $\gamma_{j-1}$ and swap the $a/b$-chain backwards along $T$, since shifting to $\gamma_{j-1}$ will change the color of $e_k$ and will hence break $T$ at an edge other than its final edge.
 
We break Case 2 into two subcases, depending on whether $T$ passes through $e_k$ or $e_0$ first. In the first subcase, we will have $T=(v_s,\ldots,v_k,u,v_0,\ldots,v_{j-1},u)$, and in the second we will have $T=(v_s,\ldots v_0,u,v_k,\ldots,v_{j-1},u)$. The second case is easy to take care of. If $T=(v_s,\ldots v_0,u,v_k,\ldots,v_{j-1},u)$, then we shift to coloring $\gamma_{j-1}$, leaving $e_{j-1}$ uncolored. Notice that $-b$ is absent at $v_k$ with respect to $\gamma_{j-1}$ (since $k >j-1$). We now perform the $a/b$-swap at $v_{j-1}$, which travels along the trail $(v_{j-1},\ldots,v_k)$ and terminates at $v_k$ (because $-b$ is absent at $v_k$). We can now color $v_{j-1}$ with the color $a$ and we are done.

Now we consider the case where $T=(v_s,\ldots,v_k,u,v_0,\ldots,v_{j-1},u)$. We note that the argument given in the previous paragraph will not work in this case. To see why, consider what happens if we shift to $\gamma_{j-1}$ and swap the $a/b$-chain at $v_{j-1}$. Since we have shifted to $\gamma_{j-1}$, the edge $uv_{k+1}$ is now colored $-b$. Thus, the swap will travel along the trail $T'=(v_{j-1},\ldots, v_0, u, v_{k+1},\ldots)$. We have not specified a second endpoint for $T'$ for good reason---we simply do not know where $T'$ will end. In fact, it is possible that the last two vertices of $T'$ are $v_j$ and $u$ so that performing the swap makes $a$ present at $u$.

Thus, to tackle the case where $T=(v_s,\ldots,v_k,u,v_0,\ldots,v_{j-1},u)$, we shall instead do this: first, shift to the coloring $\gamma_k$, leaving $e_k$ uncolored and $-b$ absent at $v_k$. Perform the $-a/-b$-swap at $v_k$, which travels along the trail $(v_k,\ldots,v_s)$. This leaves $-a$ absent at $v_k$, but this is still not quite what we want. We need to force $a$ to be absent at $v_k$ rather than $-a$. We will call the coloring that we have obtained $\gamma_k \,'$.

Next, we consider the $-a/a$-chain at $v_k$ with respect to $\gamma_k \,'$, denoted by $A$. One of two things may happen. First, $A$ may terminate at a vertex other than $u$. In this case, perform the $-a/a$-swap at at $v_k$ and extend by coloring $e_k$ with color $a$. The other possibility is that $A$ terminates at $u$. If it does, then $A=(v_k,...,v_0,u)$. Notice that when we modified $\gamma_k$ to obtain $\gamma_k \,'$, we interchanged the color present at $v_k$ from $-a$ to $-b$ and the color present at $v_s$ from $b$ to $a$, but we did not change the present or absent colors at any of the other $v_i$. Thus, with respect to $\gamma_k \,'$, we can still shift the colors on all edges below $e_k$ upwards, leaving $e_0$ uncolored. It is also important here that $e_{k-1}$ is not colored $-a$ with respect to $\gamma_k \,'$, or else shifting the colors on the edges upwards would make $-a$ present at $v_k$. Fortunately it is impossible that $e_{k-1}$ is colored $-a$, because $e_0$ is colored $-a$ and $0<j-1<k$, so that $k \geq 2$. 

Once we have shifted the colors upwards on the edges below $e_k$ with respect to $\gamma_{k}\,'$, we simply perform the $-a/a$-swap at $v_0$. This swap travels along the trail $(v_0,\ldots,v_k)$ (i.e., backwards along $A$), and terminates at $v_k$. We now color $e_0$ with the color $a$, and we are done. The part of Case 2 where $j-1 < k$ and where $T=(v_s,\ldots,v_k,u,v_0,\ldots,v_{j-1},u)$ is illustrated in Figures \ref{fig8} and \ref{fig9}.
\end{proof}
\end{lem}
We will now prove the zero-free version of signed Vizing's Theorem. We reiterate that the upper bound of $\Delta(\Sigma)+2$ appears due to the fact that a zero-free coloring always uses an even number of colors.
\begin{thm}[Zero-free Signed Vizing's Theorem]\label{vizing}Let $\Sigma$ be a signed simple graph. Then $\Delta(\Sigma) \leq \chi'(\Sigma) \leq \Delta(\Sigma)+1$ if $\Delta(\Sigma)$ is odd, and $\Delta(\Sigma) \leq \chi'(\Sigma) \leq \Delta(\Sigma)+2$ if $\Delta(\Sigma)$ is even.
\begin{proof} The proof is by induction on the number of edges, with the result being clear for signed graphs on $0$, $1$, or $2$ edges. Suppose $\Sigma$ has $n \geq 3$ edges. Choose an edge $e$ of $\Sigma$ and consider $\Sigma{\setminus}e$. If $\Delta(\Sigma) = \Delta(\Sigma{\setminus}e)$, then we obtain a coloring of $\Sigma{\setminus}e$ with either $\Delta(\Sigma)+1$ or $\Delta(\Sigma)+2$ colors (depending on the parity of $\Delta(\Sigma)$) by induction. If $\Delta(\Sigma) -1 = \Delta(\Sigma{\setminus}e)$, then we obtain a coloring of $\Sigma{\setminus}e$ with either $\Delta(\Sigma)$ or $\Delta(\Sigma)+1$ colors. Clearly a zero-free coloring with $\Delta(\Sigma)$ colors can be transformed into a coloring with $\Delta(\Sigma)+2$ colors (just add two more colors to the color set), so in either case we have a zero-free coloring of $\Sigma{\setminus}e$ such that each vertex has at least one absent color. We call this coloring $\gamma_0$.
\begin{figure}[h!]
\centering
\includegraphics[scale=.056]{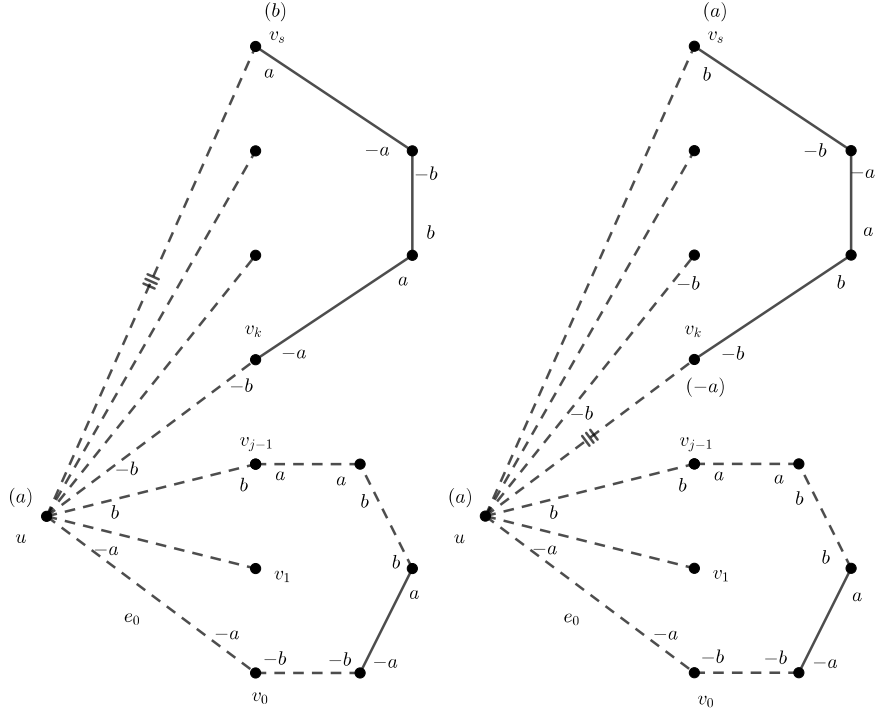}
\caption[An illustration of Case 2 in Theorem \ref{goodvizing}---part 1.]{On the left we have $\gamma_s$ along with $T=(v_s,\ldots,v_k,u,v_0,\ldots,v_{j-1},u)$, which intersects the fan at three edges, $e_k$, $e_0$, and $e_{j-1}$ (in order). On the right, we have shifted up to coloring $\gamma_k$ and then swapped the $-a/-b$-chain at $v_k$ to obtain the new coloring $\gamma_k\,'$. From the picture on the right, we proceed by attempting to swap the $-a/a$-chain at $v_k$. The details of this second swap are illustrated in Figure \ref{fig9}.  \label{fig8}}
\end{figure} 

\begin{figure}[h!]
\centering
\includegraphics[scale=.056]{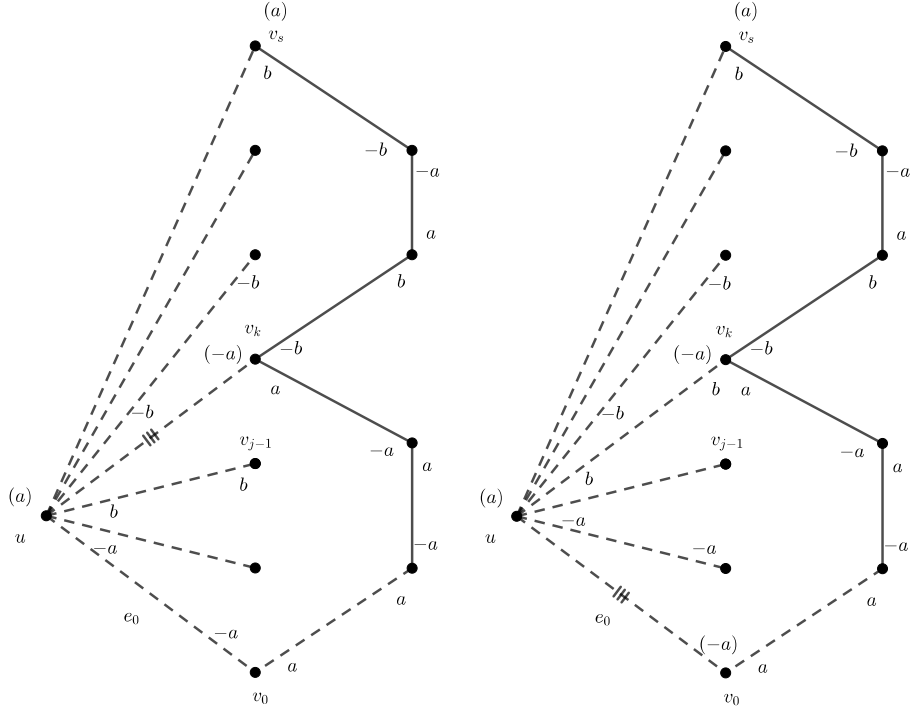}
\caption[An illustration of Case 2 in Theorem \ref{goodvizing}---part 2.]{A continuation of Figure \ref{fig8}. On the left, we see the modified coloring $\gamma_k\,'$ along with the $-a/a$-chain at $v_k$ that ends at $u$ via edge $e_0$. On the right, we have shifted the colors on the edges below $e_k$ upwards leaving $e_0$ uncolored. We are now free to swap $-a/a$ at $v_0$, which will terminate at $v_k$. This allows us to color $e_0$ with the color $a$. \label{fig9} }
\end{figure}

Start a fan $F$ at $e{:}uv_0$ using $\gamma_0$ as the initial coloring. Write $\gamma_0,\ldots,\gamma_s$ for the shifted colorings. If colors of the same magnitude are absent at both $u$ and $v_0$ with respect to $\gamma_0$, then we can apply Lemma \ref{fan} to $F$ and obtain the desired coloring of $\Sigma$. Therefore, we assume that colors of the same magnitude are not absent at $u$ and $v_0$. 

Let $b$ be absent at $v_s$ with respect to $\gamma_0$. Shift to coloring $\gamma_s$. In a similar fashion to the proof of Lemma \ref{fan}, we have an $a/b$-chain $T$ starting at $v_s$ that must end at a $b$-colored edge of $F$. Let $e_k=uv_k$ be the first edge where $T$ intersects $F$, so that either $T=(v_s,\ldots,v_k,u,\ldots,u)$, or $T=(v_s,\ldots u, v_k,\ldots, u)$, or simply $T=(v_s,\ldots,v_k, u)$. In the first two cases we are saying that $T$ passes through $e_k$ before leaving and returning to $u$ (in two different ways), and in the third case we are saying that $T$ passes through $e_k$ and then immediately ends. Notice that the third case occurs if and only if $\gamma_s(e_k)=b$.

If $\gamma_s(e_k)=b$, shift to $\gamma_k$ and perform the $a/b$-swap at $v_k$, which propagates along the trail $(v_k,...,v_s)$ and terminates at $v_s$. We can then color $e_k$ with the color $a$, finishing this case.

If $\gamma_s(e_k) = -a$, then $-a$ is absent at $v_k$ with respect to $\gamma_k$. Thus, in this case we may begin a new fan 
$F'$ with initial edge $e_k{:}uv_k$ and initial coloring $\gamma_k$. The key is that colors of the same magnitude absent at $u$ and $v_k$ with respect to $\gamma_k$. We then apply Lemma \ref{fan} to complete this case.

If $\gamma_s(e_k) = -b$, we consider two different cases. First, if $T=(v_s,\ldots,v_k,u,\ldots,u)$, then we shift to $\gamma_k$ and perform the $-a/-b$-swap at $v_k$, which travels along the trail $(v_k,\ldots v_s)$. The edge $e_k$ is now uncolored with colors of the same magnitude absent at its endpoints, and we apply Lemma \ref{fan}. 

In the second case, suppose that $T=(v_s,\ldots,u,v_k,\ldots,u)$. Since $e_k$ is colored $-b$ and $T$ passes through $u$ first and $v_k$ second, we see that no edge of $F$ is colored $-a$. If an edge of $F$ were colored $-a$, then $T$ would have to pass through it first in order to get to $e_k$, but $e_k$ is the first edge where $T$ and $F$ intersect. Thus, $T$ intersects $F$ at two edges, namely $e_k$ and $e_j$ (which is the final edge of $T$ and is colored $b$). Moreover, $T=(v_s,\ldots, u,v_k, \ldots, v_j , u)$. So, if $k >j$, we shift to $\gamma_j$ (which leaves $e_j$ uncolored and $-b$ absent at $v_k$) and perform the $a/b$-swap at $v_j$, which moves along the trail $(v_j,...,v_k)$ and terminates at $v_k$. We can then color $e_j$ with $a$ and we have finished. If $k <j$ then we shift to coloring $\gamma_k$, and perform the $-a/-b$-swap at $v_k$, which travels along $(v_k,...,v_j)$ and terminates at $v_j$. Then, we have $e_k$ uncolored with colors of the same magnitude absent at its endpoints. We apply Lemma \ref{fan} to extend the coloring. 

In any case, we are able to color all of $\Sigma$, proving the theorem.
\end{proof}
\end{thm}

We will now present another theorem that allows us to deal with the color $0$. This theorem will be pivotal in the proof of Theorem \ref{goodvizing}. Given a signed graph $\Sigma$, we write $M(\Sigma)$ for the subgraph induced by all vertices of maximum degree.
\begin{thm}\label{indeplem} Suppose $\Sigma$ is a signed graph with $\Delta(\Sigma)$ even. If $M(\Sigma)$ is an independent set, then $\Sigma$ admits an edge coloring with $\Delta(\Sigma)$ colors. 
\begin{proof} We proceed by induction on the number of vertices contained in $M(\Sigma)$. For a base case, suppose that $|M(\Sigma)|=1$. Delete an edge $e$ such that $e$ is incident with the vertex of maximum degree $u$. Then, $\Delta(\Sigma{\setminus}e)=\Delta(\Sigma)-1$ and hence there is a $\Delta(\Sigma)$-coloring of $\Sigma{\setminus}e$ by Theorem \ref{vizing}. We shall call this coloring $\gamma$.

Now, since $u$ is maximum degree, no neighbor of $u$ is maximum degree by assumption. Therefore, there is at least one absent color with respect to $\gamma$ at $u$ and at all neighbors of $u$. Thus, we are able to start a fan $F$ with initial edge $e=e_0$ and with $u$ as the hinge. From this point on, the proof is identical to that of Theorem \ref{vizing}---if colors of the same magnitude are absent at $u$ and $v_0$ then we apply Lemma \ref{fan}, and if not, we follow the process described in Theorem \ref{vizing} to extend the coloring.

 Now we proceed by induction, supposing that $|M(\Sigma)| \geq 2$ and that the statement is true for all graphs with smaller sets of maximum degree vertices. Once again, let $e$ be an edge incident with a vertex of maximum degree $u$. By induction we obtain a $\Delta(\Sigma)$-coloring of $\Sigma{\setminus}e$. We can once again build a fan with hinge $u$ and initial edge $e=e_0$ in the same way as described above, mimicing the proof of Theorem \ref{vizing}.
\end{proof}
\end{thm}

Finally, we show how the bound of $\Delta(\Sigma)+2$ can be lowered in the case where $\Delta(\Sigma)$ is even. This will complete the proof of Theorem \ref{goodvizing}.

\begin{thm} \label{vizingzero}Suppose $\Sigma$ is a signed graph with $\Delta(\Sigma)$ even. Then $\chi'(\Sigma) \leq \Delta(\Sigma)+1$. 
\begin{proof}Suppose that $\Delta(\Sigma)$ is even. Remove a maximal matching $N$ from $M(\Sigma)$. If $\Delta(\Sigma{\setminus}N) = \Delta(\Sigma)-1$, then Theorem \ref{vizing} gives a $\Delta(\Sigma)$-coloring of $\Sigma{\setminus}N$. If $\Delta(\Sigma{\setminus}N) = \Delta(\Sigma)$, then Theorem \ref{indeplem} gives a $\Delta(\Sigma)$-coloring of $\Sigma{\setminus}N$, since removing a maximal matching $N$ from $M(\Sigma)$ that does not cover all of $M(\Sigma)$ will leave $M(\Sigma{\setminus}N)$ as an independent set. 
In either case, we have obtained a zero-free coloring of $\Sigma{\setminus}N$ using $\Delta(\Sigma)$ colors. We are now free to color every edge of $N$ with the color $0$, resulting in a proper coloring of $\Sigma$ using the prescribed number of colors. 
\end{proof}

\end{thm}
In parallel with the ordinary Vizing's Theorem, the signed version of Vizing's Theorem partitions signed graphs into two classes. Indeed, we say that a signed graph $\Sigma$ is \emph{class $1$} if it admits a coloring that achieves the lower bound of $\Delta(\Sigma)$ colors. Otherwise, a signed graph is called \emph{class $2$}.

It is possible to produce two different signed graphs on the same underlying graph, one of which is class $1$ and the other class $2$---for example, positive and negative circles. However, every graph in $[\Sigma]$ has the same class as $\Sigma$. We define the \emph{class ratio} of an unsigned graph $\Gamma$ to be the number of signatures on $\Gamma$ such that the resulting signed graph is $\Delta$-colorable, divided by $2^m$, the number of possible signatures on $\Gamma$. The class ratio is denoted by $\mathcal{C}(\Gamma)$. Thus for example $\mathcal{C}(C_n)=1/2$, as only balanced circles are $2$-colorable. Since each switching class on $\Gamma$ contains the same number of signatures, $\mathcal{C}(\Gamma)$ can also be computed by counting the ratio of $\Delta$-colorable switching classes.

\subsection{Snarks}

In ordinary edge coloring, a \emph{snark} is a connected isthmus-free $3$-regular graph that does not admit an edge coloring with $3$ colors. Snarks are known for being quite hard to find, with relatively few examples known. A \emph{signed snark} is a connected isthmus-free $3$-regular signed graph that is not $3$-colorable. By Theorem \ref{goodvizing}, all signed snarks are $4$-colorable. Every unsigned snark corresponds to a signed snark on the same graph---simply sign all the edges as negative.

A natural question one might ask is the following: are the any connected isthmus-free $3$-regular graphs $\Gamma$ such that every signed graph on $\Gamma$ is a signed snark? We answer this question negatively.

\begin{lem}Let $\Gamma$ be an ordinary connected isthmus-free $3$-regular graph. There exists a signature $\sigma$ on $\Gamma$ such that $\Sigma=(\Gamma, \sigma)$ is $3$-colorable. In other words, $\mathcal{C}(\Gamma)>0$.

\begin{proof} We apply Petersen's theorem---every connected $3$-regular isthmus-free graph has a perfect matching. Let $M$ be a perfect matching in $\Gamma$. Then $\Gamma {\setminus} M$ is a $2$-regular graph---it is a union of circles. We choose $\sigma$ such that $\Gamma {\setminus}M$ is balanced, and we choose the signature of $M$ arbitrarily. Thus, we can color $\Sigma{\setminus}M$ with $\pm a$ and color $M$ with $0$, obtaining a $3$-coloring of $\Sigma$. 
\end{proof}
\end{lem}

The opposite question to that posed above is also interesting: are there any connected $3$-regular isthmus-free graphs such that $\mathcal{C}(\Gamma)=1$?

Before giving an example, we introduce a helpful concept. The \emph{frustration index} of a signed graph $\Sigma$ is the minimum number of negative edges that occur over all signed graphs in $[\Sigma]$. Equivalently, the frustration index is the minimum number of edges that must be deleted to obtain a balanced signed graph. The \emph{maximum frustration} of a graph $\Gamma$ is the maximum frustration index over all possible signatures. 

\begin{exm} Every signature of $K_{3,3}$ is $3$-colorable. In other words, $\mathcal{C}(K_{3,3})=1$.

\begin{proof} It sufficies to explain that any signature on $K_{3,3}$ contains a positive $6$-circle. The complement of the $6$-circle is a matching, so we can color the $6$-circle with $\pm a$ and the matching with $0$.

Indeed, let $\Sigma=(K_{3,3}, \sigma)$. It is known (see \cite{bowlin}) that the maximum frustration of $K_{3,3}$ is $2$. Thus, we assume that $\Sigma$ has $2$ or less negative edges. If $\Sigma$ has $0$ or $1$ negative edges, simply choose a matching that contains them. The complement of this matching is a balanced $6$-circle. If $\Sigma$ has $2$ negative edges they must be non-adjacent, since if they are adjacent we can switch to a signature with $1$ negative edge. There is a perfect matching containing any two non-adjacent edges of $K_{3,3}$.
\end{proof}
\end{exm}

In fact there is nothing particularly special about $K_{3,3}$ here---if $\Gamma$ has maximum frustration $2$ and if there is a perfect matching covering any two given edges, then every signature of $\Gamma$ is $\Delta$-colorable. Thus $\mathcal{C}(K_4)=1$, for example. 

We close this section by giving a bound on the class ratio for Hamiltonian connected $3$-regular isthmus-free graphs.

\begin{exm} If $\Gamma$ is a Hamiltonian connected $3$-regular isthmus-free graph, then $\mathcal{C}(\Gamma) \geq 1/2$. 
\begin{proof} Let $H$ be a Hamilton circle in $\Gamma$. Then the complement of $H$ is a perfect matching. We note that $H$ is positive in precisely half of the possible signatues on $\Gamma$. So, $H$ can be $2$-colored in precisely half the signatures. Thus, $\Gamma$ is $3$-colorable in at least half of its possible signatures.
\end{proof}
\end{exm}
The above argument can be generalized. If $\Gamma$ contains a $2$-regular spanning subgraph $K$ with $k$ components, then $K$ is balanced in exactly $1/2^k$ of the possible signatures. Thus in this case $\mathcal{C}(\Gamma) \geq 1/2^k$. 

\subsection{Class Ratio of Complete Graphs}

It is desirable to calculate $\mathcal{C}$ for certain classes of well known graphs. Here we briefly mention one of the simplest possibilies, $\mathcal{C}(K_n)$. We have already seen that $\mathcal{C}(K_2)=1$, $\mathcal{C}(K_3)=1/2$, and $\mathcal{C}(K_4)=1$. Based on this evidence one might guess that $\mathcal{C}(K_n) = 1$ if $n$ is even, and $\mathcal{C}(K_n)=1/2$ if $n$ is odd. However, in reality the situation is not quite this simple.

First, we note that if $n$ is odd then $\mathcal{C}(K_n)\leq1/2$. This is because any $\Delta$-coloring of $K_n$ (where $n$ is odd) is a decomposition of $K_n$ into $\Delta/2$ balanced $2$-regular spanning subgraphs. If the signature has an odd number of negative edges, then some circle in the decomposition must contain an odd number of negative edges and hence be unbalanced. Exactly half of all possible signatures have an odd number of negative edges.

In order to prove that $\mathcal{C}(K_n)\geq1/2$ when $n$ is odd it would suffice to show that any signature with an even number of negative edges can be $\Delta$-colored. However, this is not true---consider $K_5$ with an all-negative signature. This signature has $10$ negative edges, but a decomposition of $K_5$ into two balanced spanning subgraphs must be a decomposition into two circles of length $5$. Neither of these circles will be positive as they each contain $5$ negative edges. Thus, it is in fact true that $\mathcal{C}(K_5) < 1/2$. We leave it as an open problem to determine a precise formula for $\mathcal{C}(K_n)$.

\section{Line Graphs}

In this section we show that every signed edge coloring can be realized as a vertex coloring of a signed line graph. This is a desirable property for signed edge coloring to possess, since unsigned edge coloring posseses the very same property. We recall that a \emph{vertex coloring} of $\Gamma$ is an assignment of a color to each of the vertices of $\Gamma$. Such a coloring is \emph{proper} if no two adjacent vertices have the same color. 

To vertex color a signed graph, we assign a color from the set $M_n$ to each of its vertices. We employ the definition of propriety discovered by Zaslavsky \cite{zaslav2}---a \emph{proper signed vertex coloring} has the requirement that positive edges do not have the same color at their endpoints, and negative edges do not have colors with the same magnitude and opposite sign at their endpoints. 

One of the nice features of this definition is that it extends to switching classes in a natural way. If $c$ is a proper vertex coloring of $\Sigma$ and $\Sigma \sim \Sigma'$ via switching vertex set $X$, then we can obtain a proper coloring $c'$ of $\Sigma'$ by simply negating $c(x)$ for all $x \in X$. In this way $c$ generates a proper vertex coloring for each member of $[\Sigma]$.

\subsection{Bidirected Graphs}

The easiest way to define the line graph of a signed graph is through the use of bidirected graphs. A \emph{bidirected graph} is a pair $(\Gamma, \tau)$, where $\Gamma$ is a graph and $\tau:I(\Gamma)\rightarrow \{+,-\}$ is a \emph{bidirection}. When $\tau(v,e)=+$ we imagine an arrow drawn on $e$ that points into $v$, and when $\tau(v,e)=-$ we imagine an arrow drawn on $e$ that points away from $v$. An edge $e$ is \emph{extraverted} if both of its $\tau$ values are $+$, \emph{introverted} if both of its $\tau$ values are $-$, and \emph{coherent} otherwise.

A negation of the $\tau$ values for a certain edge $e$ is a \emph{reorientation} of $e$. Thus a reorientation of an extraverted edge is an introverted edge and vice versa, while a reorientation of a coherent edge remains coherent. Reorientation is an equivalence relation on bidirected graphs, and hence there is a partition of the set of bidirected graphs into \emph{reorientation classes}. The reorientation class of a given bidirected graph $B$ is denoted by $\vec{B}$.

Bidirected graphs can be thought of as orientations of signed graphs. For a given bidirected graph $B=(\Gamma, \tau)$, there is a natural associated signed graph $\Sigma_B = (\Gamma, \sigma_{\tau})$, obtained by setting $\sigma_{\tau}(e) = -\tau(v,e)\tau(w,e)$ for $e{:}{vw}$. We say that $B$ is an \emph{orientation} of $\Sigma_B$. In other words, positive edges correspond to coherent edges, and negative edges correspond to extraverted and introverted edges.

\begin{figure}[h!]
\centering
\includegraphics[scale=.064]{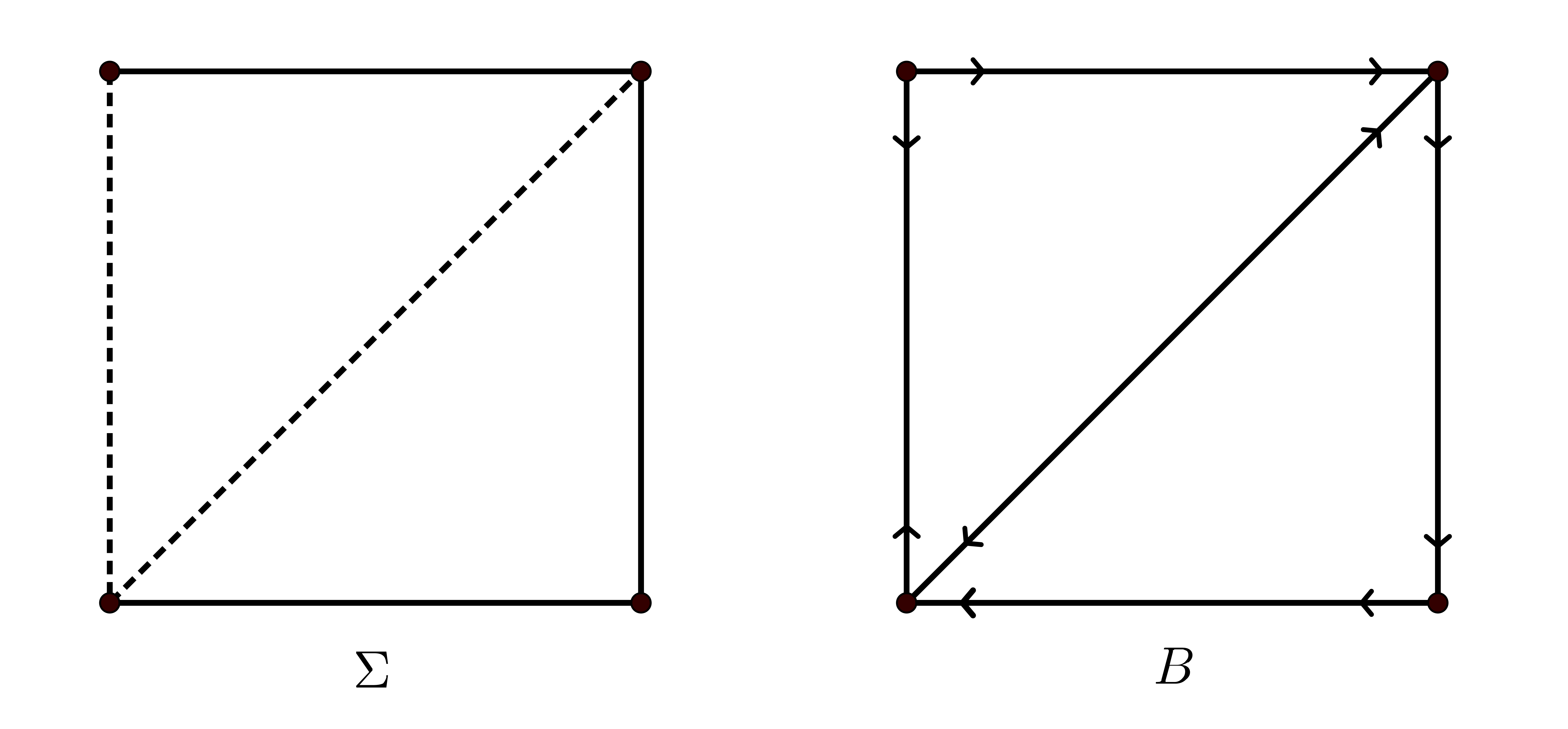}
\caption[A signed graph $\Sigma$ and one of its orientations $B$.]{On the left, a signed graph $\Sigma$. On the right, $B$, one of its $32$ possible orientations. \label{fig10} }
\end{figure}

A signed graph with $m$ edges has $2^m$ possible orientations, and each of these orientations can be obtained from any other by reorientation of the appropriate edges. Thus, if $\Sigma$ has a single orientation $B$, then $\vec{B}$ contains exactly the $2^m$ orientations of $\Sigma$. See Figure \ref{fig10} for a picture of a signed graph and an orientation.

In light of the fact that bidirected graphs are orientations of signed graphs, we may use the terminology of signed graphs to refer to bidirected graphs when it is not confusing to do so. For example, a subgraph of a bidirected graph is \emph{balanced} if its correspoding signed subgraph is balanced. 

A vertex $v$ of a bidirected graph can be \emph{switched} by negating all $\tau$ values involving $v$. It is easy to check that bidirected switching is compatible with signed switching as far as orientations go---switching both $B$ and $\Sigma_B$ at $v$ will result in $B'$ and $\Sigma_B'$ such that $B'$ is an orientation of $\Sigma_B'$.

\subsection{Coloring Bidirected Graphs}

The definition of edge coloring for a signed graph cooperates nicely with bidirected graphs. We define an edge coloring of a bidirected graph in the following way. Notice that our definition is in terms of edges rather than in terms of incidences.

\begin{defn}An $n$-edge coloring (or more briefly, an $n$-coloring) $\gamma$ of a bidirected graph $B$ is an assignment of colors from $M_n$ to each edge of $B$. Such a coloring is \emph{proper} if $\tau(v,e)\gamma(e)\neq \tau(v,f) \gamma(f)$ for all edges $e$ and $f$ that are adjacent at vertex $v$. 
\end{defn}

Thus in a proper coloring if the $\tau$ values are equal the edges may not have the same color, while if they are not equal the edges may not have opposite colors. 

The purpose of this definition is to enable us to view an edge coloring of $\Sigma$ as an edge coloring of one of its orientations.

\begin{lem}\label{orientcolor}Suppose $\Sigma$ is a signed graph and $B$ is one of its orientations. Suppose $\gamma$ is an edge coloring of $\Sigma$. Then there is a unique edge coloring $\gamma_B$ of $B$ such that $\tau(v,e)\gamma_B(e)=\gamma(v,e)$ for all $v$ and incident $e$. 

\begin{proof} We describe how to define $\gamma_B$ for each edge of $B$. Take any edge $e$ and let $v$ be one of its endpoints. We set $\gamma_B(e)=\tau(v,e)\gamma(v,e)$. Thus, $\tau(v,e)\gamma_B(e)=\gamma(v,e)$. We now must check the other endpoint of $e$. Let $w$ be the other endpoint of $e$. Then $\tau(w,e)\gamma_B(e) = -\sigma(e)\tau(v,e)\gamma_B(e)=-\sigma(e)\gamma(v,e)=\gamma(w,e)$. 
\end{proof}
\end{lem}

Furthermore, if we are given an edge coloring $\gamma_B$ of $B$, we can uniquely recover the edge coloring $\gamma$ of $\Sigma$ by setting $\gamma(v,e)=\tau(v,e)\gamma_B(e)$. 

The following lemma shows that a coloring of one orientation automatically generates a coloring of every other possible orientation.

\begin{lem}\label{orientcolor2} Suppose $\gamma_B$ is an edge coloring of $B$. Let $B'$ be a reorientation of $B$. Then there is a unique edge coloring $\gamma_{B'}$ such that $\gamma_{B'}(e) \tau'(v,e) = \gamma_B(e) \tau(v,e)$ for all incidences $(v,e)$. 

\begin{proof} Suppose $e$ is an edge of $B$, and suppose that $e$ is one of the edges that gets reoriented when passing from $B$ to $B'$. In this case we define $\gamma_{B'}(e):= -\gamma_B(e)$. If $e$ is not reoriented, we set $\gamma_{B'}(e) = \gamma_B(e)$. Either way, $\gamma_{B'}(e) \tau'(v,e) = \gamma_B(e) \tau(v,e)$.
\end{proof}
\end{lem}

Suppose $\Sigma$ is a signed graph with two different orientations $B$ and $B'$. Let $\gamma$ be an edge coloring of $\Sigma$, and let $\gamma_{B}$ and $\gamma_{B'}$ be the corresponding edge colorings of $B$ and $B'$ in the sense of Lemma \ref{orientcolor}. Then, Lemma \ref{orientcolor2} tells us that we can obtain $\gamma_{B}$ from $\gamma_{B'}$ by negating the colors on the edges that must be reoriented to change from $B$ to $B'$. 

Thus, $\gamma$ uniquely determines $\gamma_{B}$ for every orientation $B$ of $\Sigma$, and conversely $\gamma_{B}$ determines all other $\gamma_{B'}$, each of which determine $\gamma$. 

The following figure illustrates Lemma \ref{orientcolor} and Lemma \ref{orientcolor2}.
\begin{center}
\begin{figure}[h!]\label{fig11}
\centering
\includegraphics[scale=.056]{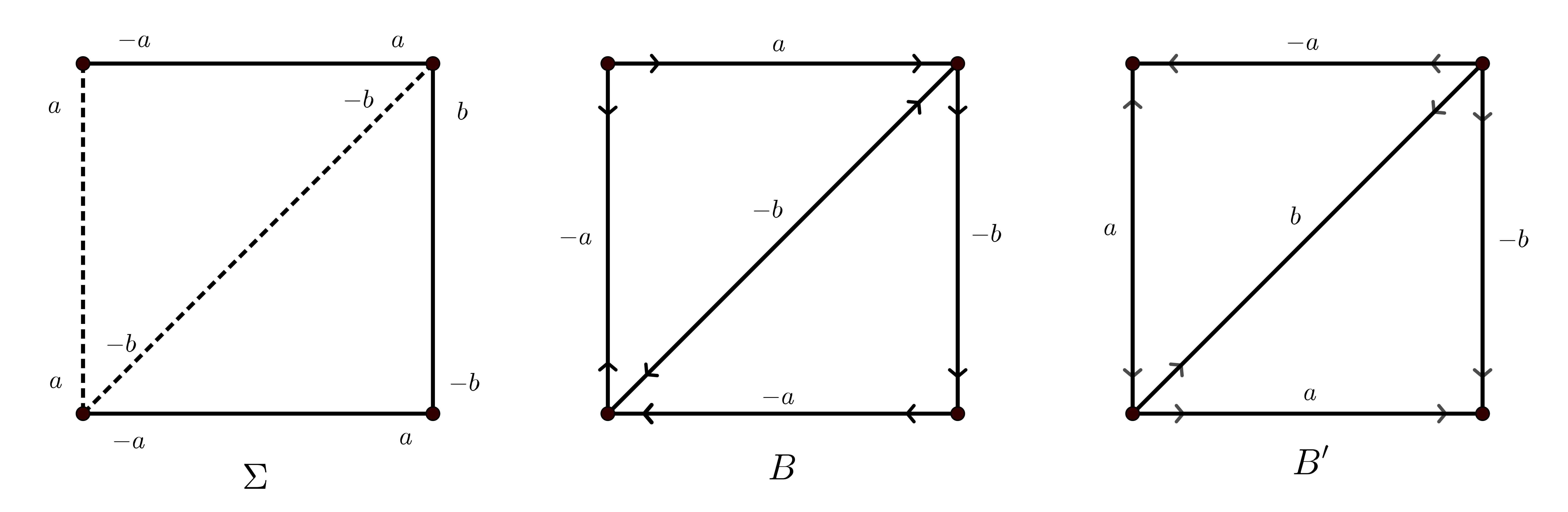}
\caption[$\Sigma$, two of its orientations, and two edge colorings.]{From left to right; $\Sigma$ along with an edge coloring, one of its orientations $B$ with the corresponding bidirected edge coloring, and a second orientation $B'$ with its corresponding edge coloring. We can pass from $B$ to $B'$ by negating colors on the edges that have been reoriented.}
\end{figure} 
\end{center}

\subsection{Defining Line Graphs}

First, we present the definition of the signed line graph, as originally studied by Zaslavsky \cite{zaslav4}. Recall that the \emph{line graph} of an ordinary graph $\Gamma$ is the graph $L(\Gamma)$ whose vertices are the edges of $\Gamma$, and two vertices of $L(\Gamma)$ are adjacent in $L(\Gamma)$ if and only the corresponding edges are adjacent in $\Gamma$. Thus, $L(\Gamma)$ is the graph of edge adjacency for $\Gamma$. If $e$ is an edge of $\Gamma$, we write $\ell_e$ to represent the corresponding vertex in $L(\Gamma)$. 

\begin{defn} The \emph{line graph} of a bidirected graph $B=(\Gamma,\tau)$ is $L(B):=(L(|B|),\vec{\tau})$, where $\vec{\tau}(\ell_e,\ell_e\ell_f):=\tau(v,e)$ (where $v$ is the common vertex of edges $e$ and $f$ in $B$). 
\end{defn}

Thus for example the line graph of an all-extraverted $B$ is itself all extraverted. The main purpose of the bidirected line graph is to act as a tool that enables us to define the line graph of a signed graph. Before proceeding with the definition, we point out that the line graph of $\Sigma$ turns out to be a switching class of signed graphs, rather than a single signed graph. 

\begin{defn}[Signed Line Graph]\label{slg} The \emph{line graph} of a signed graph $\Sigma$ is obtained by the following procedure:
\begin{enumerate}
\item Choose any orientation $B$ of $\Sigma$.
\item Find $L(B)$, the bidirected line graph of $B$.
\item Find the signed graph corresponding to $L(B)$. Denote this by $\Sigma_{L(B)}$. 
\item The line graph of $\Sigma$ is defined as the switching class of $\Sigma_{L(B)}$. We write $\Lambda(\Sigma) = [\Sigma_{L(B)}]$.
\end{enumerate}
\end{defn}

\begin{center}
\begin{figure}[h!]
\centering
\includegraphics[scale=.063]{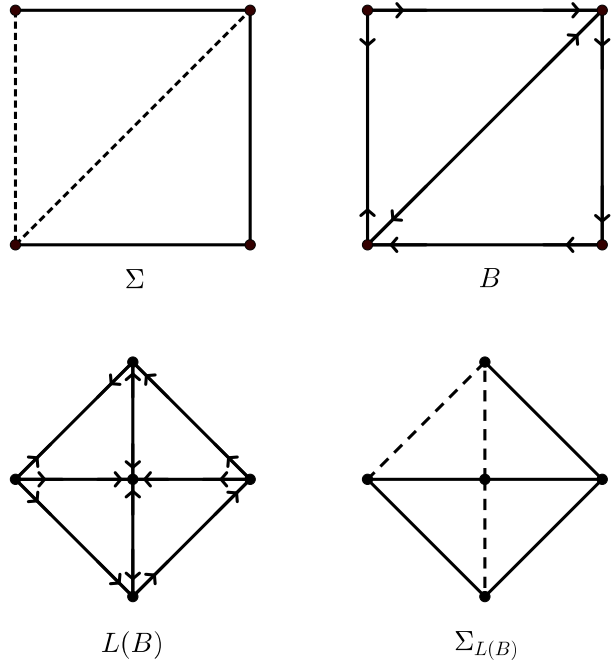}
\caption[A signed graph $\Sigma$, one of its orientations $B$, the bidirected line graph $L(B)$ of $B$, and the signed graph corresponding to $L(B)$, $\Sigma_{L(B)}$.]{A signed graph $\Sigma$, one of its orientations $B$, the bidirected line graph $L(B)$ of $B$, and the signed graph corresponding to $L(B)$, $\Sigma_{L(B)}$. The line graph of $\Sigma$ is $\Lambda(\Sigma)=[\Sigma_{L(B)}]$, the switching class of $\Sigma_{L(B)}$. A reorientation of one of the edges of $B$ will switch the corresponding vertex in $L(B)$, which switches the corresponding vertex in $\Sigma_{L(B)}$.\label{fig12}}
\end{figure} 
\end{center}

It is important to notice that the above definition does not depend on the choice of $B$. Indeed, reorienting an edge $e$ of $B$ will have the effect of switching the vertex $\ell_e$ in $L(B)$. Since $\Lambda(\Sigma)$ is a switching class, reorienting $e$ does not change the line graph. An illustration of Definition \ref{slg} is given in Figure \ref{fig12}.

Not only is the line graph of $\Sigma$ a switching class, but every signed graph switching equivalent to $\Sigma$ has the same line graph as $\Sigma$.

\begin{lem} If $\Sigma \sim \Sigma'$, then $\Sigma$ and $\Sigma'$ have the same line graph.

\begin{proof} We switch a single vertex $v$ of $\Sigma$ and observe what effect it has on the line graph. Let $e$ and $f$ be two edges that are adjacent at $v$. Indeed, switching $v$ has the effect of negating all $\tau$ values at $v$ for any orientation $B$ of $\Sigma$. Thus, in $L(B)$, the edge $\ell_e\ell_f$ gets reoriented. This reorientation has no effect on $\Lambda(\Sigma)$, and hence does not change the line graph. 
\end{proof}
\end{lem}

\subsection{Edge Coloring in Terms of Vertex Coloring}

We will now study how one can interpret an edge coloring of $\Sigma$ in terms of a vertex coloring of the line graph of $\Sigma$. It turns out that an edge coloring of $\Sigma$ corresponds with a vertex coloring (in Zaslavsky's sense) of $[-\Sigma_{L(B)}]$, the negative of the line graph of $\Sigma$. 

\begin{thm}\label{lineg} There is a bijection between (proper) edge colorings of $\Sigma$ and (proper) vertex colorings (in Zaslavsky's sense) of $[-\Sigma_{L(B)}]=-\Lambda(\Sigma)$, the negative of the line graph of $\Sigma$.

\begin{proof} Let $\gamma$ be a (proper) edge coloring of $\Sigma$. We will describe how to define a (proper) vertex coloring $c$ of $-\Lambda(\Sigma)$ in terms of $\gamma$.

Choose any orientation $B$ of $\Sigma$, using $\gamma$ to induce the unique edge coloring $\gamma_B$ of $B$. Since the edges of $B$ are the vertices of $L(B)$, we can think of $\gamma_B$ as a vertex coloring of $L(B)$. In turn, we think of $\gamma_B$ as a vertex coloring of $-\Lambda(\Sigma)$. Thus, we have a bijection between (not necessarily proper) edge colorings of $\Sigma$ and vertex colorings of $-\Lambda(\Sigma)$.

We now wish to prove that $\gamma$ is a proper edge coloring of $\Sigma$ if and only if $\gamma_B$ is a proper vertex coloring of $-\Lambda(\Sigma)$. Indeed, let edges $e$ and $f$ be adjacent at vertex $v$ in $\Sigma$. Since $\gamma$ is proper, $\gamma(v,e)\neq \gamma(v,f)$. Equivalently, in $B$, $\gamma_B(e) \tau(v,e) \neq \gamma_B(f) \tau(v,f)$. Thus, equivalently in $L(B)$, the previous rule becomes $\gamma_B(\ell_e) \tau(\ell_e, \ell_e\ell_f) \neq \gamma_B(\ell_f) \tau(\ell_f, \ell_e\ell_f)$. Thus, when passing to $-\Lambda(\Sigma)$, we have $\gamma_B(\ell_e) \neq \sigma(\ell_e \ell_f) \gamma_B(\ell_f)$, which is precisely the definition of a proper vertex coloring. 
\end{proof}
\end{thm}

\begin{center}
\begin{figure}[h!]
\centering
\includegraphics[scale=.063]{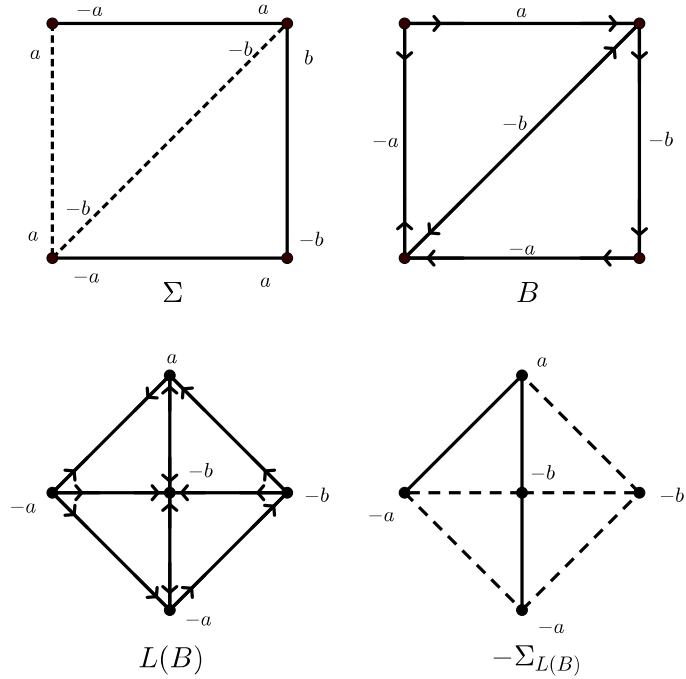}
\caption[A proper edge coloring of $\Sigma$ and the corresponding vertex coloring of its line graph.]{A proper edge coloring of $\Sigma$, the corresponding proper edge coloring of $B$, the corresponding vertex coloring of $L(B)$, and the corresponding proper vertex coloring of $-\Sigma_{L(B)}$. \label{fig13}}
\end{figure} 
\end{center}

Figure \ref{fig13} shows a signed graph along with an edge coloring and the corresponding vertex coloring of the negative of its line graph. The intermediate bidirected graph steps are shown as well.

Theorem \ref{lineg} specializes to ordinary graphs nicely when $\Sigma$ is all-negative. Indeed, if $\Sigma$ is all-negative we can choose $B$ so that it is all-extraverted. Then $L(B)$ is also all-extraverted, so that $-\Sigma_{L(B)}$ is all-positive. Thus, edge colorings of an all-negative $\Sigma$ correspond to vertex colorings of an all-positive $-\Sigma_{L(B)}$. This is what we expect, since edge colorings of an all-negative $\Sigma$ correspond to ordinary edge colorings, and vertex colorings of an all-positive $-\Sigma_{L(B)}$ correspond to ordinary vertex colorings.

\section{Additional Topics}
\subsection{Reversibility and The Linear Arboricity Conjecture}
The \emph{linear arboricity} $\la(\Gamma)$ of an ordinary simple graph is the minimum number of linear forests (acyclic subgraphs of maximum degree $2$) into which its edges can be partitioned. The linear arboricity of $\Gamma$ is at least $\lceil \Delta / 2 \rceil$ since each linear forest in such a partition uses at most two edges incident with a vertex of maximum degree. In 1981, Akiyama et al.\ conjectured \cite{lineararb} that the linear arboricity of a graph is bounded above by $\lceil (\Delta+1) / 2 \rceil$, a conjecture that remains unresolved at present. Since $\lceil \Delta/2 \rceil$ and $\lceil (\Delta +1) /2 \rceil$ are either equal or consecutive integers, the linear arboricity conjecture states that $\la(\Gamma)$ is either $\lceil \Delta/2 \rceil$ or $\lceil (\Delta +1 )/2 \rceil$---a dichotomous statement that is reminiscent of Vizing's Theorem. As we are about to see, the problem of computing Linear Arboricity can be naturally phrased in terms of signed graph edge coloring.

Suppose $\Sigma$ is a signed graph, and $\gamma$ is a proper edge coloring. We say that an edge $e$ is \emph{reversible with respect to $\gamma$} if $e$ lies in a path component of $\Sigma_{\gamma(e)}$. If $e$ is reversible with respect to $\gamma$, we can negate the sign of $e$ and easily find another proper coloring $\gamma \,'$ such that the magnitudes of the colors in $\gamma$ are the same as those in $\gamma \,'$. This occurs because we can always color any bidirected path with only two colors. 

We say that the proper coloring $\gamma$ is \emph{completely reversible} if every edge of $\Sigma$ is reversible with respect to $\gamma$. The following two lemmas should now be evident.

\begin{lem} Let $\Sigma$ be a signed graph. A proper edge coloring $\gamma$ is completely reversible if and only if every component of every magnitude graph is a path. 
\end{lem}

\begin{lem} If $\Sigma=(\Gamma,\sigma)$ admits a completely reversible $n$-coloring, then so does every signed graph with underlying graph $\Gamma$.

\end{lem}

Thus, in a completely reversible coloring each $\Sigma_a$ is a linear forest, and conversely any partition of the edges of $\Sigma$ into linear forests yields a completely reversible coloring. We note that if $a=0$ this linear forest is a matching (not an arbitrary linear forest), so we look only at zero-free colorings. We write $\chi'_R(\Sigma)$ for the minimum number of colors needed in a completely reversible zero-free proper coloring of $\Sigma$. We can now phrase the linear arboricity conjecture in terms of edge coloring.

\begin{cnj}[Linear Arboricity Conjecture] For any simple signed graph $\Sigma$, $\Delta(\Sigma) \leq \chi'_R(\Sigma) \leq \Delta(\Sigma) +2$. 
\end{cnj}

This is equivalent to the linear arboricity conjecture since each $\Sigma_a$ requires two colors  but is one linear forest.

\subsection{A Change of Convention}

If we tweak our definition of signed edge coloring slightly we end up with $\Sigma_a$ graphs that differ from those we see normally. In fact, the change that we make to the definition will turn the $\Sigma_a$ into arbitrary antibalanced subgraphs, which links edge coloring to a problem that has already been studied---balanced decomposition. 

The change is this: when we define a proper edge coloring, we insist that $\gamma(v,e) \neq -\gamma(v,f)$ for all edges $e$ and $f$ adjacent at vertex $v$. This is the negative of the usual definition, so let us call such a coloring \emph{antiproper}. In the language of orientations, an antiproper coloring is one that satisfies $\gamma(e) \tau(v,e) \neq - \gamma(f) \tau(v,f)$. 

We first note that the Vizing-style lower bound of $\Delta(\Sigma)$ does not hold for an antiproper coloring, since every incidence at a given vertex may be colored the same. However, despite this, the magnitude graphs $\Sigma_a$ still do have some nice structure.

\begin{lem} Let $\gamma$ be an antiproper edge coloring of $\Sigma$. Then each $\Sigma_a$ graph is antibalanced. Furthermore, any antibalanced graph can be antiproperly colored with $2$ colors ($\pm a$). 

\begin{proof} Let $C$ be a circle in $\Sigma_a$ and suppose $C$ has $m$ edges. Switch $\Sigma$ so that $C$ consists of a negative path of length $m-1$, and one additional edge $e$ that is either positive or negative depending on the sign of $C$. Since $\gamma$ is antiproper, every edge in the path $C{\setminus}e$ must have both its incidences colored (without loss of generality) $a$. Thus, $e$ must be negative as well, or else we do not have an antiproper coloring. Hence, $C$ switches to all negative and so $\Sigma_a$ is antibalanced. 

To prove the other statement, let $A$ be an arbitrary antibalanced signed graph. Switch $A$ to all negative and color every incidence as $a$. This gives a $2$-coloring (not a $1$-coloring, since we must use $M_2=\{-a,a\}$). 
\end{proof}
\end{lem}

Thus, an antiproper coloring of $\Sigma$ is equivalent to a decomposition of $\Sigma$ into antibalanced subgraphs. These antibalanced subgraphs are arbitrary (except in the case of $\Sigma_0$, which must be a matching). In turn, a decomposition of $\Sigma$ into antibalanced subgraphs corresponds to a decomposition of $-\Sigma$ into balanced subgraphs. 

In \cite{zaslav1}, Zaslavsky studies the \emph{balanced decomposition number} of $\Sigma$---the smallest number of balanced sets into which its edges can be partitioned, denoted by $\delta_0(\Sigma)$. The balanced decomposition number is a parameter that encapsulates how far a given signed graph is from being balanced---the higher the balanced decomposition number, the ``less balanced'' $\Sigma$ is. 

The balanced decomposition number is a generalization of the \emph{biparticity} $\beta_0$ of an unsigned graph---the fewest number of bipartite sets into which the edges can be partitioned. In particular, $\beta_0(\Gamma) = \delta_0(-\Gamma)$ (here $-\Gamma$ means an all negatively signed $\Gamma$), since the balanced subgraphs in an all negative graph are precisely the bipartite subgraphs. Biparticity is known to be connected to the chromatic number $\chi(\Gamma)$ by the formula 
$$\beta_0(\Gamma) = \lceil \log_2(\chi(\Gamma)) \rceil,$$
discovered independently by Harary-Hsu-Miller \cite{bipar1} and Matula \cite{bipar2}.
A similar theorem for $\delta_0$ was given by Zaslavsky in \cite{zaslav1}.
\begin{thm}[Zaslavsky's Balanced Decomposition Theorem] \label{zaslav2} If $\Sigma$ has at least one edge, $\delta_0(\Sigma)=\lceil \log_2( \chi^*(-\Sigma))\rceil$, where $\chi^*(\Sigma)$ is the zero-free vertex chromatic number of $\Sigma$.

\end{thm}
Let us write $\chi_A'$ for the minimum number of colors needed in any zero-free antiproper coloring. Then Zaslavsky's Theorem immediately gives $\chi_A'(\Sigma) = 2\delta_0(-\Sigma)=2\lceil \log_2( \chi^*(\Sigma))\rceil$. The reason for the multiplication by $2$ is the fact that in any antibalanced decomposition, each antibalanced set requires two colors.

We close this section by offering an interesting interpretation of antiproper colorings in terms of a line graph. In contrast with proper colorings, we do not have to negate the line graph to obtain the correspondence. The proof is evident from our discussion of signed line graphs. 

\begin{thm} Antiproper edge colorings of $\Sigma$ correspond to proper vertex colorings of the line graph $\Lambda(\Sigma)$.
\end{thm}

\subsection{Total Coloring} In this section we will discuss how to define total coloring for a signed graph. There are a couple of ways that we can do this, both of which are interesting in their own right. 

Recall that a \emph{total coloring} of an ordinary graph $\Gamma$ is an assignment of colors to its vertices and edges such that no two adjacent vertices, adjacent edges, or incident vertices and edges share a color. In other words, a total coloring is simultaneously a proper vertex coloring and a proper edge coloring, and the interaction between them is that incident vertices and edges also do not share a color.

We would like to make a similar definition for signed graphs, and ideally, our definition should specialize to the ordinary definition in some way. However, there is a subtlety that we must deal with---signed vertex colorings correspond to ordinary vertex colorings when $\Sigma$ is balanced, but signed edge colorings correspond to ordinary edge colorings when $\Sigma$ is antibalanced. Thus, we are compelled to make the following definition:

\begin{defn} A \emph{signed total coloring} $\mu$ of $\Sigma$ is an assignment of colors from $M_n$ to the vertices and incidences of $\Sigma$ such that:
\begin{enumerate}
\item $\mu$ restricted to $V(\Sigma)$ is a proper vertex coloring of $-\Sigma$.
\item $\mu$ restricted to $I(\Sigma)$ is a proper edge coloring. 
\item $\mu(v) \neq \mu(v,e)$ for all incident vertices and edges $v$ and $e$.
\end{enumerate}
\end{defn}

Thus, a signed total coloring of an all-negative $\Sigma$ corresponds to a total coloring of $|\Sigma|$, which is what we desire. Since signed vertex and edge coloring are both compatible with switching, signed total coloring is too.

\begin{lem} Let $\mu$ be a signed total coloring of $\Sigma$. Let $\Sigma' \sim \Sigma$ via switching function $\eta$. Then $\mu'$ is a total coloring of $\Sigma'$, where $\mu'$ is obtained from $\mu$ by negating the colors on all vertices and incidences that were switched via $\eta$.

\begin{proof} This follows immediately from the compatibility of vertex and edge coloring with switching.
\end{proof}
\end{lem}

Let $\chi''(\Sigma)$ be the total chromatic number---the fewest number of colors used in any signed total coloring of $\Sigma$. It is easy to notice that $\chi''(\Sigma) \geq \Delta(\Sigma) +1$, as a maximum degree vertex requires $\Delta$ different colors on its incident edges and one additional color for itself. 

It is natural to look for a Vizing-style upper bound for $\chi''$, and indeed for ordinary graphs it has been long conjectured that the upper bound is $\Delta +2$, although no proof has been found. This is known as the \emph{total coloring conjecture}, first posed by Behzad. Interestingly, many people attribute this conjecture to Vizing, but according to Shahmohamad \cite{shah}, Behzad is indeed the sole author of the conjecture. The upper-bound of $\Delta+2$ is known to hold for some specific classes of ordinary graphs, such as $r$-partite graphs \cite{yap}, $3$-regular graphs \cite{rosen}, and most planar graphs. 

Based on our results concerning signed Vizing's theorem, we conjecture the following. 

\begin{cnj}[Signed Total Coloring Conjecture] Any simple signed graph $\Sigma$ admits a total coloring with $\Delta(\Sigma)+2$ colors.
\end{cnj}

It is left as an open problem to prove various special cases of this conjecture, such as $r$-partite graphs and $3$-regular graphs.

Let us revisit the definition of total coloring. In the definition, we insisted on having a proper vertex coloring of $-\Sigma$ so that total coloring would specialize for antibalanced graphs. Instead, let us do the following:

\begin{defn} A \emph{twisted signed total coloring} $\mu$ of $\Sigma$ is an assignment of colors from $M_n$ to the vertices and incidences of $\Sigma$ such that:
\begin{enumerate}
\item $\mu$ restricted to $V(\Sigma)$ is a proper vertex coloring of $\Sigma$.
\item $\mu$ restricted to $I(\Sigma)$ is a proper edge coloring. 
\item $\mu(v) \neq \mu(v,e)$ for all incident vertices and edges $v$ and $e$.
\end{enumerate}
\end{defn}

The difference in definitions is only a single negative sign, but something interesting happens. If $\Sigma$ is all-negative, a twisted total coloring corresponds to an ordinary edge coloring of $|\Sigma|$. However, if $\Sigma$ is all-positive, a twisted total coloring corresponds to an ordinary vertex coloring of $|\Sigma|$. Thus, a twisted total coloring specializes in two different ways, but not necessarily at the same time. In fact, both ways may happen at the same time.

\begin{lem} A twisted total coloring of $\Sigma$ corresponds to a total coloring of $|\Sigma|$ if and only if $\Sigma$ is both balanced and antibalanced (i.e., if and only if $\Sigma$ is balanced and bipartite). 
\begin{proof} We require $\Sigma$ to switch to both all-positive and all-negative.
\end{proof}
\end{lem}


\begin{thebibliography}{99}
\bibitem{lineararb} J. Akiyama, G. Exoo, and F. Harary, \emph{Covering and packing in graphs. IV. Linear Arboricity}, Networks, \textbf{11}(1) (1981), 69--72.
\bibitem{bowlin} G. Bowlin, \emph{Maximum Frustration in Bipartite Signed Graphs}, Elec. L. of Comb. \textbf{19}(4) (2012), 1--13.
\bibitem{harary} F. Harary, \emph{On the notion of balance of a signed graph}, Michigan Math. J. 2 (1953), no. 2, 143--146.
\bibitem{bipar1} F. Harary, D. Hsu and Z. Miller, \emph{The biparticity of a graph}, J. Graph Theory (1977), 1: 131--133. 
\bibitem{maca} E. M\'a\v{c}ajov\'a, A Raspaud, and M. \v{S}koviera, \emph{The chromatic number of a signed graph}, Elec. J. of Comb. \textbf{23}(1) (2016), 1--10. 
\bibitem{bipar2} D. Matula, \emph{$k$-Components, clusters and slicings in graphs}, SIAM J. Appl. Math. \textbf{22}(3) (1972), 459--480. 
\bibitem{rosen} M. Rosenfeld, \emph{On the total coloring of certain graphs}, Israel J. Math. \textbf{9}(3), 396--402.
\bibitem{shah} H. Shahmohamad, \emph{The history of the total chromatic number conjecture}, https://arxiv.org/pdf/1104.3170.
\bibitem{vz} V. Vizing, \emph{On an Estimate of the Chromatic Class of a $p$-graph}, Diskret. Analiz \textbf{3} (1964), 23--30.
\bibitem{yap} H. Yap, \emph{Total Colourings of Graphs}, Bull. London Math. Soc. \textbf{21}(2) (1989), 159--163.
\bibitem{zaslav1} T. Zaslavsky, \emph{Balanced decompositions of a signed graph}, J. Combin. Theory Ser. B \textbf{43}(1) (1987), 1--13.
\bibitem{zaslav2} T. Zaslavsky, \emph{Signed graph coloring}, Discrete Math. \textbf{39}(2) (1982), 215--228.
\bibitem{zaslav3} T. Zaslavsky, \emph{Signed graphs}, Discrete Appl. Math. 4 (1982), 47--74.
\bibitem{zaslav4} T. Zaslavsky, \emph{Matrices in the theory of signed simple graphs}, Advances in Discrete Mathematics and Applications, (Ramanujan Math. Soc. Lect. Notes Mysore, India), \textbf{13} (2010), 207--229.





\end{thebibliography}
\end{document}